\documentclass[11pt]{amsart}

\usepackage[margin=2.5cm]{geometry}
\usepackage[ruled,vlined]{algorithm2e}
\usepackage{booktabs}
\usepackage{graphicx}
\usepackage[T1]{fontenc}
\usepackage[utf8]{inputenc}
\usepackage{url}
\usepackage{txfonts}
\usepackage{listings}
\usepackage{hyperref}
\usepackage{tikz}
\usepackage{tikz-qtree}
\usepackage{subfig}
\lstdefinelanguage{GAP}{%
    morekeywords=[2]{and,break,continue,do,elif,else,end,fail,false,fi,for,%
        function,if,in,local,mod,not,od,or,rec,repeat,return,then,true,%
        until,while, %
        AffineSemigroup, Difference, List, Intersection, CatenaryDegreeOfElementInNumericalSemigroup, TameDegreeOfElementInNumericalSemigroup, NumericalSemigroup,IsEmpty, Filtered, ForAny, Gaps, Generators, Rank, Minimum},%
    sensitive=true,%
    morecomment=[l]\#,%
    morestring=[b]',%
    morestring=[b]",%
    basicstyle=\ttfamily,
    }%

\newtheorem{theorem}{Theorem}[section]
\newtheorem{proposition}[theorem]{Proposition}
\newtheorem{lemma}[theorem]{Lemma}

\newtheorem{corollary}[theorem]{Corollary}

\theoremstyle{definition}
\newtheorem{definition}[theorem]{Definition}
\newtheorem{remark}[theorem]{Remark}
\newtheorem{example}[theorem]{Example}

\begin{document}

\title[Generalized Numerical semigroups up to isomorphism]{Generalized Numerical semigroups up to isomorphism}

%\date{}
\author{CARMELO CISTO}
\address{Universit\`{a} di Messina, Dipartimento di Scienze Matematiche e Informatiche, Scienze Fisiche e Scienze della Terra, Viale Ferdinando Stagno D'Alcontres 31, 98166 Messina, Italy}
\email{carmelo.cisto@unime.it}

\author{GIOIA FAILLA}
\address{Universit\`a Mediterranea di Reggio Calabria, DIIES\\
Via Graziella, Feo di Vito,  Reggio Calabria, Italy}
\email{gioia.failla@unirc.it}

\author{FRANCESCO NAVARRA}
	\address{Sabanci University, Faculty of Engineering and Natural Sciences, Orta Mahalle, Tuzla 34956, Istanbul, Turkey}
\email{francesco.navarra@sabanciuniv.edu}

\thanks{ }



\begin{abstract}
A generalized numerical semigroup is a submonoid $S$ of $\mathbb{N}^d$ with finite complement in it. We characterize isomorphisms between these monoids in terms of permutation of coordinates. Considering the equivalence relation that identifies the monoids obtained by the action of a permutation and establishing a criterion to select a representative from each equivalence class, we define some procedures for generating the set of all generalized numerical semigroups of given genus up to isomorphism. Finally, we present computational data and explore properties related to the number of generalized numerical semigroups of a given genus up to isomorphism.
\end{abstract}

%, effective generator, genus, monomial order, algorithm

%\maketitle
\maketitle

\section{Introduction}

Let $\mathbb{N}$ be the set of non negative integers. A \emph{ numerical semigroup}  is a submonoid $S$ of $\mathbb{N}$ having finite complement in it. This notion is widely studied in literature, also with connections in different fields of mathematics, as commutative algebra, algebraic geometry or coding theory.  For a good overview we refer to \cite{ns-app} and \cite{rosales2009numerical}. In \cite{failla2016algorithms}, a straightforward generalization of numerical semigroups is provided for submonoids in $\mathbb{N}^d$, with $d$ a positive integer. In particular, the authors introduce the concept of \emph{generalized numerical semigroup} (GNS), that is a submonoid $S$ of $\mathbb{N}^d$ such that the set $\operatorname{H}(S)=\mathbb{N}^d \setminus S$ is finite. The set $\operatorname{H}(S)$ is called the set of \emph{gaps} (or \emph{holes}) of $S$ and its cardinality $\operatorname{g}(S)=\vert\operatorname{H}(S)\vert$ is called the \emph{genus} of $S$. Generalized numerical semigroups and some generalizations of them are studied in several recent papers (see for instance \cite{bernardini_atoms, bernardini_corner, bhardwaj1, bhardwaj2, can2024irreducible, mahir_wilf,cisto2022some}). In \cite{failla2016algorithms} Failla et al. develop an algorithm, inspired by \cite{bras2008fibonacci, bras2009bounds}, that leads to the construction of a tree, that we denote by $\mathcal{T}_{d,\preceq}$, having as set of vertices the set $S_{d}$ of  all  generalized  numerical semigroups in $\mathbb{N}^d$. In particular, if $g$ is a positive integer and $\mathcal{S}_{g,d}$ is the set of GNSs in $\mathbb{N}^d$ having genus $g$, the construction allows to produce all semigroups in $\mathcal{S}_{g+1,d}$ from all semigroups in $\mathcal{S}_{g,d}$. Moreover, the behavior of the values $n_{g,d}=|S_{g,d}|$ is analyzed for low values of $d$ and $g$. A first implementation of this procedure in a software is considered in \cite{garcia2018extension}. New algorithms to deal with GNSs are introduced in \cite{cisto2021algorithms}. To be precise, it is shown how to compute the set $\mathcal{S}_{g,d}$, the set of gaps from a minimal generating set of a GNS and the set of minimal set of generators of a GNS, once its gaps are known. In particular, the set $\mathcal{S}_{g,d}$ is computed by using a different procedure of \cite{failla2016algorithms}, and it is performed by the construction of a different tree, that we denote by $\mathcal{T}'_{g,d,\preceq}$, having $\mathcal{S}_{g,d}$ as set of vertices. We mention that the structure of the trees $\mathcal{T}_{d,\preceq}$ and $\mathcal{T}'_{g,d,\preceq}$ depends on the definition of a relaxed monomial order $\preceq$. By implementing  these algorithms in the computer algebra system \texttt{GAP} \cite{GAP} and with the package \texttt{numericalsgps} \cite{numericalsgps}, some values of $n_{g,d}$ are computed, for instance for $d=2$ up to genus $g=21$, for $d=3$ up to genus $g=13$ and for higher dimesion up to genus $g=9$.

In this paper we introduce the notion of isomorphism between generalized numerical semigroups. We consider the equivalence relation $\simeq$ in the set $\mathcal{S}_{g,d}$, where if $S_1,S_2\in \mathcal{S}_{g,d}$, then $S_1 \simeq S_2$ if $S_2$ can be obtained by a permutation of the coordinates of $S_1$. We notice that if $S_1 \simeq S_2$, then $S_1$ and $S_2$ are isomorphic as monoids. We show that the converse holds, that is, if $S_1$ and $S_2$ are isomorphic GNSs in $\mathbb{N}^d$, then $S_1 \simeq S_2$. 
We first consider the set of GNSs $S\subseteq \mathbb{N}^d$ such that $|[S]_\simeq|=1$, that we call \emph{equivariant}, showing how it is possible to arrange the set of equivariant GNSs in $\mathbb{N}^d$ in a tree. Fixed a positive integer $k$, this construction allows to produce the set of all equivariant GNSs having genus $g\leq k$. Furthermore, considering the tree $\mathcal{T}_{d,\preceq}$, a natural question is how it is possible to ``trim'' this tree in order to produce only one semigroup in $[S]_\simeq$, for every $S\in \mathcal{S}_{d}$, and such that the trimmed tree preserves the same behaviour of the first one in its construction. The same question arises for the tree $\mathcal{T}'_{g,d,\preceq}$. To reach this aim, we determine a possible criterion to choose a representative in the equivalence class $[S]_\simeq$, for every $S\in \mathcal{S}_{g,d}$. Using this criterion, we show that it is possible to use the same technique of the construction of $\mathcal{T}_{d,\preceq}$, in order to generate only the representatives of all semigroups in $\mathcal{S}_{d}$. We also investigate if the same choice of representatives allows to generate all the representatives of the semigroups in $\mathcal{S}_{g,d}$, using the construction of the tree $\mathcal{T}'_{g,d,\preceq}$, obtaining that it is possible only for some choices of the monomial order $\preceq$.
 Finally, let $N_{g,d}$ be the cardinality of the set of all generalized numerical semigroups in $\mathcal{S}_{g,d}$ up to isomorphism.  By implementations of the  procedures mentioned above, we provide some values of $N_{g,d}$ and we prove that $N_{g,g}=N_{g,d}$ for all $d\geq g$. 

The structure of the paper is as follows. Section 2 consists of basic definitions, results and techniques for generalized numerical semigroups.
Section 3 is devoted to prove that two GNSs in $\mathbb{N}^d$ are isomorphic if and only if they differ only by a permutation of coordinates. This allow us to state that, in the set of all GNS in $\mathbb{N}^d$, the equivalence relation, where all GNSs obtained by a permutation of coordinates are identified, corresponds to the isomorphism relation. In Section 4 we introduce the concept of equivariant GNS, that is a GNS such that its isomorphism equivalence class has cardinality equal to one. We give some properties and a procedure to generate this type of GNSs. In Section 5 we study some procedure to generate the set of all GNSs in $\mathbb{N}^d$ with a fixed genus up to isomorphism. To reach this aim, we establish a possible criterion to identify a representative in every isomorphism equivalence class of a GNS. Finally, in Section 6, we provide some properties about the sequence $N_{g,d}$ of the number of GNSs in $\mathbb{N}^d$ with genus $g$ and give some computational data of it, for different values of $d$ and $g$.

\section{Preliminaries and notations}

A \emph{generalized numerical semigroup} (GNS) is a submonoid $S$ of $\mathbb{N}^d$ with finite complement in $\mathbb{N}^d$. The set $\operatorname{H}(S)$ is called the set of \emph{gaps} (or \emph{holes}) of $S$ and its cardinality $\operatorname{g}(S)=|\operatorname{H}(S)|$ is called the \emph{genus} of $S$. This notion has been introduced in \cite{failla2016algorithms} as a  generalization of the well known concept of numerical semigroup. By known results, every submonoid $S$ of $\mathbb{N}^d$ admits a unique minimal set of generators $\operatorname{A}(S)$, that is, every element in $S$ is a linear combination of element in $\operatorname{A}(S)$ with coefficient in $\mathbb{N}$ and no proper subset of $\operatorname{A}(S)$ has this property. In particular, $A(S)=S^{*}\setminus(S^{*}+S^{*})$, where $S^{*}=S\setminus \{0\}$. It is known that if $S$ is a GNS, then $\operatorname{A}(S)$ is finite and we denote $\operatorname{e}(S)=|\operatorname{A}(S)$|, called embedding dimension of $S$.

\noindent Recall that the natural partial order $\leq$ in $\mathbb{N}^d$ is defined by $\mathbf{a}\leq \mathbf{b}$ if $\mathbf{b}-\mathbf{a}\in \mathbb{N}^d$. As usual, when we define an order relation $\preceq$ on a set $A$, which is a total order, and we consider $x,y \in A$, the notation $x\prec y$ means that $x\preceq y$ and $x\neq y$. In the following definition we provide a class of total orders in $\mathbb{N}^d$, introduced in \cite{failla2016algorithms}, that will be useful in this work.

\begin{definition}
A \emph{relaxed monomial order} is a total order $\preceq$ in $\mathbb{N}^d$ satisfying the two properties:
\begin{itemize}
\item [a)] $\mathbf{0}\preceq \mathbf{v}$ for all $\mathbf{v}\in \mathbb{N}^d$;
\item[b)] for all $\mathbf{v},\mathbf{w}\in \mathbb{N}^d$ such that $\mathbf{v}\prec \mathbf{w}$ then $\mathbf{v}\prec \mathbf{w}+\mathbf{u}$ for all $\mathbf{u}\in \mathbb{N}^d$.
\end{itemize}
The total order $\preceq$ is called a \emph{monomial order} ff the second property is replaced with
\begin{itemize}
\item [b')]  for all $\mathbf{v},\mathbf{w}\in \mathbb{N}^d$ such that $\mathbf{v}\prec \mathbf{w}$ then $\mathbf{v}+\mathbf{u}\prec \mathbf{w}+\mathbf{u}$ for all $\mathbf{u}\in \mathbb{N}^d$,
    \end{itemize}
\end{definition}    

% \noindent Recall that a \emph{relaxed monomial order} (introduced in \cite{failla2016algorithms}) is a total order $\preceq$ in $\mathbb{N}^d$ satisfying the two properties: $a)$ $\mathbf{0}\prec \mathbf{v}$ for all $\mathbf{v}\in \mathbb{N}^d$; $b)$ for all $\mathbf{v},\mathbf{w}\in \mathbb{N}^d$ such that $\mathbf{v}\prec \mathbf{w}$ then $\mathbf{v}\prec \mathbf{w}+\mathbf{u}$ for all $\mathbf{u}\in \mathbb{N}^d$. If the second property is replaced with $b')$ for all $\mathbf{v},\mathbf{w}\in \mathbb{N}^d$ such that $\mathbf{v}\prec \mathbf{w}$ then $\mathbf{v}+\mathbf{u}\prec \mathbf{w}+\mathbf{u}$ for all $\mathbf{u}\in \mathbb{N}^d$, $\prec$ is called a \emph{monomial order}.

\begin{example}\label{def:orders}
For any $\mathbf{a} = (a_1, a_2, \ldots, a_d)$ and $\mathbf{b} = (b_1, b_2, \ldots, b_d)$ in $\mathbb{N}^d$, consider the following orders:
\begin{itemize}
\item[a.] The \emph{Lexicographic Order} is defined by $\mathbf{a} \prec_{\text{lex}} \mathbf{b}$ if and only if there exists an index $i \in \{1, 2, \ldots, d\}$ such that $a_j = b_j$ for all $j < i$ and $a_i < b_i$.

\item[b.] The \emph{Graded Lexicographic Order} is defined by $\mathbf{a} \prec_{\text{glex}} \mathbf{b}$ if and only if one of the following occurs:
\begin{enumerate}
    \item $\sum_{j=1}^{d} a_j < \sum_{j=1}^{d} b_j$.
    \item $\sum_{j=1}^{d} a_j = \sum_{j=1}^{d} b_j$ and $a \prec_{\text{lex}} b$. 
\end{enumerate}
%\item[c.] Let $\prec$ a monomial order in $\mathbb{N}^d$. Define $a \prec_{\text{min}} b$ if and only if one of the following occurs:
%\begin{enumerate}
%    \item $\min\{a_1,\ldots,a_d\} < \min\{b_1,\ldots,b_d\}$.
%    \item $\min\{a_1,\ldots,a_d\} < \min\{b_1,\ldots,b_d\}$ and $a \prec b$.
%    \end{enumerate}
\item[c.] Let $\preceq$ be a monomial order in $\mathbb{N}^d$ and $\leq$ be the natural partial order in $\mathbb{N}^d$. Denote $X_d=\{(x_1,\ldots,x_d) \in \mathbb{N}^d\mid x_i=0\text{ for all }i\in \{1,\ldots,d-1\}\}$. Define $\mathbf{a} \prec_{\text{1}} \mathbf{b}$ if and only if one of the following occurs:
\begin{enumerate}
    \item $\mathbf{a},\mathbf{b}\in X_d$ and $\mathbf{a}\leq \mathbf{b}$.
    \item $\mathbf{a}\in X_d$ and $\mathbf{b}\notin X_d$.
    \item $\mathbf{a},\mathbf{b}\notin X_d$ and $\mathbf{a}\prec \mathbf{b}$.
    \end{enumerate}  
    \end{itemize}
The previous total orders are example of relaxed monomial orders in $\mathbb{N}^d$. Moreover, the lexicographic order and the graded lexicographic order are monomial orders, while the third is not a monomial order in general. In fact, considering $\preceq$ equal to Graded Lexicographic Order, in $\mathbb{N}^3$ we have $(0,0,3)\prec_1 (0,1,1)$ but $(0,0,3)+(1,0,0)\nprec_1 (0,1,1)+(1,0,0)$.
\end{example}

%INSERIRE UN REMARK: NON ESISTONO ELEMENTI NON NULLI DI UN GNS MINORI (RISPETTO AD UN ORDINE MONOMIALE) DI TUTTI I BUCHI.

Related to relaxed monomial orders in $\mathbb{N}^d$, we now introduce the following notions.

\begin{definition}
Let $S$ be a submonoid of $\mathbb{N}^d$ and $\preceq$ be a relaxed monomial order. We define
\begin{enumerate}
  \item $\mathbf{m}_{\preceq}(S)$ is the smallest element of $S^{*}$ with respect to $\preceq$  and it is  called \emph{multiplicity} of $S$ with respect to $\preceq$.
%  \item $e(S)$ the cardinality of the minimal set of generators $\operatorname{A}(S)$.
  \item $\mathbf{F}_{\preceq}(S)=\max_\preceq (\operatorname{H}(S)$, called the \textit{Frobenius element} of $S$ with respect to $\preceq$ (convetionally, $\mathbf{F}_{\preceq}(\mathbb{N}^d)=(-1,\ldots,-1)$).
  \item $\mathbf{U}_{\preceq}(S)=\{\mathbf{x}\in \operatorname{A}(S)\mid \mathbf{F}_{\preceq}(S)\prec \mathbf{x}\}$. %$ the set of all minimal generators greater than  Frobenius number $\mathbf{F}_{\preceq}(S)$ , with respect to $\preceq$.
\end{enumerate}

\end{definition}

We briefly recall some useful notions on oriented graphs, that we will use in this work. An oriented (or direct) graph is a structure $G=(\mathcal{V},\mathcal{E})$ such that $\mathcal{V}$ is a nonempty set, called the set of \emph{vertices}, and $\mathcal{E}\subseteq \mathcal{V}\times \mathcal{V}$, called the set of \emph{edges}. If $\mathbf{x},\mathbf{y}\in \mathcal{V}$, a \emph{path} from $\mathbf{x}$ to $\mathbf{y}$ of length $n$ is a sequence of distinct vertices $\mathbf{x}_0,\mathbf{x}_2,\ldots,\mathbf{x}_n$ such that $\mathbf{x}_0=\mathbf{x}$, $\mathbf{x}_n=\mathbf{y}$ and $(\mathbf{x}_{i-1},\mathbf{x}_i)\in \mathcal{E}$ for all $i\in\{1,\ldots,n\}$. An oriented graph is called a \emph{tree} if there exists a vertex $\mathbf{v}$, called \emph{root}, such that for every vertex $\mathbf{x}$ there exists a unique path from $\mathbf{x}$ to $\mathbf{v}$. In this case, the length of the unique path from $\mathbf{x}$ to $\mathbf{v}$ is called the \emph{depth} of $\mathbf{x}$ and if $(\mathbf{x}, \mathbf{y})$ is an edge, we say that $\mathbf{x}$ is a \emph{child} of $\mathbf{y}$. A vertex in a tree is a \emph{leaf} if it has no children. In the context of numerical semigroups and GNSs, some procedures are related to some constructions of tree oriented graphs. 

\noindent For any positive integers $d,g$, let us denote:
\begin{itemize}
    \item $\mathcal{S}_d=\{S\subseteq \mathbb{N}^d\mid S\text{ is a GNS }\}$.
    \item $\mathcal{S}_{g,d}=\{S\in \mathcal{S}_d\mid \operatorname{g}(S)=g\}$
\end{itemize}
In \cite{failla2016algorithms}, the computation of the set $\mathcal{S}_{g,d}$ is performed from the the construction of an oriented graph $\mathcal{T}_{d,\preceq}$, defined for any fixed relaxed monomial order $\preceq$. In particular, $\mathcal{T}_{d,\preceq}$ is a tree having the set of all GNSs in $\mathbb{N}^d$ as set of vertices and whose root is $\mathbb{N}^d$. Moreover, if $S$ is a GNS, then the children of $S$ are the semigroups $S\setminus \{\mathbf{g}_1\}, \ldots, S\setminus \{\mathbf{g}_m\}$, such that $\mathbf{U}_{\preceq}(S)=\{\mathbf{g}_1,\ldots, \mathbf{g}_m\}$. In this way, all GNSs of genus $g$ are generated, without redundancies, from all GNSs of genus $g-1$. Then the set of all GNSs of genus $g$ consists of all vertices in $\mathcal{T}_{d,\preceq}$ of depth $g$. We remark that, going back in a branch of the tree, every GNS $S\varsubsetneq \mathbb{N}^d$ is a child of $S\cup \{\mathbf{F}_\preceq(S')\}$ and the pairs $(S, S\cup \{\mathbf{F}_\preceq(S')\})$ is an edge of the graph. %In particular, if $\mathcal{S}_d$ is the set of all GNSs in $\mathbb{N}^d$, the building of this tree can be related to the transform

\medskip
\noindent Finally, if $S\subseteq \mathbb{N}^d$ is a GNS, the set of \emph{pseudo Frobenius} elements of $S$ is defined as the set $\operatorname{PF}(S)=\{\mathbf{h}\in \operatorname{H}(S)\mid \mathbf{h}+\mathbf{s}\in S \text{ for all }\mathbf{s}\in S^*\}$, while the set of \emph{special gaps} of $S$ is defined as $\operatorname{SG}(S)=\{\mathbf{h}\in \operatorname{PF}(S)\mid 2\mathbf{h}\in S\}$. Observe that $\operatorname{SG}(S)\neq \emptyset$, since each maximal gap with respect to the natural partial order in $\mathbb{N}^d$ is a special gap. By the definition, it is not difficult to see that the set of special gaps of a GNS characterizes the unitary extension of a GNS in this sense:

\begin{proposition}
    Let $S\subseteq \mathbb{N}^d$ be a GNS and $\mathbf{h}\in \operatorname{H}(S)$. Then $S\cup \{\mathbf{h}\}$ is a GNS if and only if $\mathbf{h}\in \operatorname{SG}(S)$. 
\end{proposition}

Finally, in the whole paper we denote by $\{\mathbf{e}_1,\ldots,\mathbf{e}_2,\ldots,\mathbf{e}_d\}$ the set of standard basis vectors of the vector space $\mathbb{R}^d$. In particular, we denote $\mathbf{e}_1=(0,\ldots,0,1)$, $\mathbf{e}_2=(0,\ldots,0,1,0)$, \ldots, $\mathbf{e}_d=(1,0,\ldots,0)$. That is, we assume that $\mathbf{e}_1\prec_{\text{lex}} \mathbf{e}_2 \prec_{\text{lex}} \cdots \prec_{\text{lex}} \mathbf{e}_d$, with reference to Example~\ref{def:orders}.

\section{Permutations and isomorphisms between GNSs}\label{iso}
If $S,T$ are monoids, we recall that a map $f:S\rightarrow T$ is an homomorphism (of monoid) if for all $\mathbf{a},\mathbf{b} \in S$ then $f(\mathbf{a}+\mathbf{b})=f(\mathbf{a})+f(\mathbf{b})$. Moreover $f$ is an \emph{isomorphism} if it is also a bijective map. Obviously if $f$ is an isomorphism then it is invertible and $f^{-1}:T\rightarrow S$ is also an isomorphism. As usual, if $f:S\rightarrow S$ is an isomorphism of $S$ in itself then $f$ is also called an \emph{automorphism} of $S$.  \\

\begin{remark}
Let $S,T$ be submonoids of $\mathbb{N}^d$ and $f:S\rightarrow T$ be an isomorphism. Then 
\begin{enumerate}
\item  $f(\mathbf{0})=\mathbf{0}$, in fact for all $\mathbf{s}\in S$ we have $f(\mathbf{s})= f(\mathbf{s}+\mathbf{0})=f(\mathbf{s})+f(\mathbf{0})$ and the claim is obtained since $T$ is cancellative.

\item If $S$ and $T$ are finitely generated then it is easy to verify that if $\mathbf{s}$ is a minimal generator of $S$ then $f(\mathbf{s})$ is a minimal generator of $T$. In particular if $S$ and $T$ are GNSs then $\operatorname{e}(S)=\operatorname{e}(T)$.
\end{enumerate}

\end{remark}

After introducing some notations, we define a particular class of isomorphisms of submonoids of $\mathbb{N}^d$. 

\begin{definition}
	Let $d\in \mathbb{N}$. We denote with $\operatorname{P}_d$ the set of permutations on $\{1,\ldots,d\}$. Let $\{\mathbf{e}_1,\ldots,\mathbf{e}_d\}$ be the set of the standard basis of the vector space $\mathbb{R}^d$. If $\mathbf{x}=\sum_i^d x_i \mathbf{e}_i\in \mathbb{N}^d$, we define $\sigma(\mathbf{x})=\sum_i^d x_i \mathbf{e}_{\sigma(i)}$. Furthermore, if $A\subseteq \mathbb{N}^d$ we denote $\sigma(A)=\{\sigma(\mathbf{a})\mid \mathbf{a}\in A\}$ and if $\mathbf{x}\in \mathbb{N}^d$, define $\mathrm{orb}(\mathbf{x})=\{\sigma(\mathbf{x})\in\mathbb{N}^d \mid \sigma \in \operatorname{P}_d\}$.
\end{definition}

It is easy to verify that if $S$ is a submonoid of $\mathbb{N}^d$, then $\sigma(S)$ is also a submonoid of $\mathbb{N}^d$ for every $\sigma \in \operatorname{P}_d$. It is easy also to prove the following result.

\begin{proposition}
Let $S$ be a submonoid of $\mathbb{N}^d$ and $\sigma\in \operatorname{P}_d$. Let $f_\sigma: S \rightarrow \sigma(S)$ be the map defined by $f_\sigma(\mathbf{s})=\sigma(\mathbf{s})$. Then $f_\sigma$ is an isomorphism.	
\end{proposition}

In this section we study the isomorphisms between GNSs in $\mathbb{N}^d$ with fixed $d$. In particular we want to show that if $S$ and $T$ are both submonoids of $\mathbb{N}^d$ with finite complement in it and they are isomorphic, then there exists $\sigma \in \operatorname{P}_d$ such that $T=\sigma(S)$. We start providing the following trivial property.

\begin{proposition}\label{prop:iso-gaps}
Let $S\subseteq \mathbb{N}^d$ be a GNS and $\sigma\in \operatorname{P}_d$. Then $\operatorname{H}(\sigma(S))=\sigma(\operatorname{H}(S))$.
\end{proposition}
\begin{proof}
It trivially follows, since for all elements $\mathbf{x}\in \mathbb{N}^d$ we have $\mathbf{x}\in S$ if and only if $\sigma(\mathbf{x})\in \sigma(S)$.
\end{proof}

If $L\subseteq \mathbb{N}^d$, the partial order $\leq_L$ in $\mathbb{N}^d$ is defined by $\mathbf{a}\leq_L \mathbf{b}$ if $\mathbf{b}-\mathbf{a}\in L$. Recall that if $S\subseteq \mathbb{N}^d$ is a GNS and $\mathbf{n}\in S\setminus \{\mathbf{0}\}$ the Ap\'ery set is defined as $\operatorname{Ap}(S,\mathbf{n})=\{\mathbf{x}\in S\mid \mathbf{x}-\mathbf{n}\notin S\}$. Moreover $\operatorname{PF}(S)=\{\mathbf{w}-\mathbf{n}\mid \mathbf{w} \in \mathrm{Maximals}_{\leq_S}\operatorname{Ap}(S,\mathbf{n})\}$ (see \cite[Proposition 1.4]{cisto2019irreducible}). 
So for each $\mathbf{x}\in \operatorname{PF}(S)$ then $\mathbf{x}=\mathbf{s}-\mathbf{n}$ with $\mathbf{n}\in S\setminus \{\mathbf{0}\}$ and $\mathbf{s}\in \mathrm{Maximals}_{\leq_S}\operatorname{Ap}(S,\mathbf{n})\}$.

\begin{lemma} \label{apery1}
Let $S,T$ be GNSs in $\mathbb{N}^d$ and $f:S\rightarrow T$ be an isomorphism. Let $\mathbf{n}\in S\setminus \{\mathbf{0}\}$, then the following holds:
\begin{enumerate}
\item $\mathbf{x}\in \operatorname{Ap}(S,\mathbf{n})$ if and only if $f(\mathbf{x})\in \operatorname{Ap}(T,f(\mathbf{n}))$.
\item $\mathbf{x}\in \mathrm{Maximals}_{\leq_S}\operatorname{Ap}(S,\mathbf{n})$ if and only if $f(\mathbf{x})\in \mathrm{Maximals}_{\leq_T}\operatorname{Ap}(T,f(\mathbf{n}))$.
\item $|\operatorname{PF}(S)|=|\operatorname{PF}(T)|$.
\end{enumerate}
\end{lemma}
\begin{proof}
(1) Let $\mathbf{x}\in S$ such that $\mathbf{x}-\mathbf{n}\notin S$. If $f(\mathbf{x})-f(\mathbf{n}) \in T$  then $f(\mathbf{x})=f(\mathbf{n})+\mathbf{t}$ with $\mathbf{t}\in T$. Since $f$ is an isomorphism, $\mathbf{t}=f(\mathbf{s})$ with $\mathbf{s}\in S$, so $f(\mathbf{x})=f(\mathbf{n})+f(\mathbf{s})=f(\mathbf{n}+\mathbf{s})$. It follows $\mathbf{x}=\mathbf{n}+\mathbf{s}$ that is $\mathbf{x}-\mathbf{n}\in S$, a contradiction. This argument shows in particular that $\mathbf{x}-\mathbf{n}\in S$ if and only if $f(\mathbf{x})-f(\mathbf{n}) \in T$, from which we obtain our claim.\\
(2) Let $\mathbf{x}\in \mathrm{Maximals}_{\leq_S}\operatorname{Ap}(S,\mathbf{n})$ and suppose that there exists $\mathbf{t} \in \operatorname{Ap}(T,f(\mathbf{n}))$ such that $f(\mathbf{x})\leq_{T} \mathbf{t}$. By the previous point $\mathbf{t}=f(\mathbf{s})$ with $\mathbf{s}\in \operatorname{Ap}(S,\mathbf{n})$ and by  $f(\mathbf{x})\leq_{T} \mathbf{t}$ and the fact that $f$ is an isomorphism we obtain $f(\mathbf{s})=f(\mathbf{x})+f(\mathbf{y})$ with $\mathbf{y}\in S$. Then $\mathbf{s}=\mathbf{x}+\mathbf{y}$, in particular $\mathbf{x}\leq_{S} \mathbf{s}$, a contradiction. The converse holds with the same argument.\\
(3) It easily follows by (2) and \cite[Proposition 1.4]{cisto2019irreducible}.
\end{proof}

\begin{lemma}
Let $S,T$ be GNSs in $\mathbb{N}^d$, $f:S\rightarrow T$ be an isomorphism and $\mathbf{x}\in \operatorname{SG}(S)$. Let $\mathbf{n}\in S\setminus \{\mathbf{0}\}$ and $\mathbf{s}\in \mathrm{Maximals}_{\leq_S}\operatorname{Ap}(S,\mathbf{n})$ such that $\mathbf{x}=\mathbf{s}-\mathbf{n}$ and consider $\mathbf{y}=f(\mathbf{s})-f(\mathbf{n})$. Then $\mathbf{y}\in \operatorname{SG}(T)$ and $f(2\mathbf{x})=2\mathbf{y}$.
\label{SG}
\end{lemma}
\begin{proof}
Let $\mathbf{x}\in \operatorname{SG}(S)$ and express $\mathbf{x}=\mathbf{s}-\mathbf{n}$ with $\mathbf{n}\in S\setminus \{\mathbf{0}\}$ and $\mathbf{s}\in \mathrm{Maximals}_{\leq_S}\operatorname{Ap}(S,\mathbf{n})$. By Lemma~\ref{apery1} and \cite[Propsition 1.4]{cisto2019irreducible} 
we have that $\mathbf{y}=f(\mathbf{s})-f(\mathbf{n})\in \operatorname{PF}(T)$. In order to prove that $\mathbf{y}\in \operatorname{SG}(T)$ it suffices to show that $2\mathbf{y}\in T$. By hypothesis $2\mathbf{x}=2\mathbf{s}-2\mathbf{n}\in S$, so $f(2\mathbf{x})+f(2\mathbf{n})=f(2\mathbf{s})$ and in particular $f(2\mathbf{x})=f(2\mathbf{s})-f(2\mathbf{n})$. By the property of $f$ and since $\mathbf{n},\mathbf{s}\in S$, we obtain $f(2\mathbf{x})=2(f(\mathbf{s})-f(\mathbf{n}))$. Therefore $f(2\mathbf{x})=2\mathbf{y}$ and since $f$ is an isomorphism we have $2\mathbf{y}\in T$, that is, $\mathbf{y}\in \operatorname{SG}(T)$.   
\end{proof}

\begin{lemma}
Let $S,T$ be GNSs in $\mathbb{N}^d$, $f:S\rightarrow T$ be an isomorphism and $\mathbf{x}\in \operatorname{SG}(S)$. Let $\mathbf{n}\in S\setminus \{\mathbf{0}\}$ and $\mathbf{s}\in \mathrm{Maximals}_{\leq_S}\operatorname{Ap}(S,\mathbf{n})$ such that $\mathbf{x}=\mathbf{s}-\mathbf{n}$ and consider $\mathbf{y}=f(\mathbf{s})-f(\mathbf{n})$. Define the map $\overline{f}:S\cup \{\mathbf{x}\}\rightarrow T\cup \{\mathbf{y}\}$ with $\overline{f}(\mathbf{x})=\mathbf{y}$ and $\overline{f}(\mathbf{s})=f(\mathbf{s})$ for all $\mathbf{s}\in S$. Then $\overline{f}$ is an isomorphism of GNSs.
\label{extend1}
\end{lemma}
\begin{proof}
Trivially $\overline{f}$ is bijective, it suffices to prove that it is an homomorphism. Let $\mathbf{x}_1,\mathbf{x}_2\in S\cup \{\mathbf{x}\}$. If $\mathbf{x}_1,\mathbf{x}_2\in S$ it is trivial $\overline{f}(\mathbf{x}_1+\mathbf{x}_2)=f(\mathbf{x}_1+\mathbf{x}_2)=f(\mathbf{x}_1)+f(\mathbf{x}_2)=\overline{f}(\mathbf{x}_1)+\overline{f}(\mathbf{x}_2)$. If $\mathbf{x}_1=\mathbf{x}_2=\mathbf{x}$, considering Lemma~\ref{SG}, we obtain $\overline{f}(2\mathbf{x})=f(2\mathbf{x})=2\mathbf{y}=\overline{f}(\mathbf{x})+\overline{f}(\mathbf{x})$. So, suppose $\mathbf{x}_1=\mathbf{x}$ and $\mathbf{x}_2\in S\setminus \{\mathbf{0}\}$ and consider that $2\mathbf{x},2\mathbf{x}_2, \mathbf{x}+\mathbf{x}_2\in S$. Therefore $\overline{f}(2\mathbf{x}+2\mathbf{x}_2)=f(2\mathbf{x})+f(2\mathbf{x}_2)=2\mathbf{y}+2f(\mathbf{x}_2)=2(\overline{f}(\mathbf{x})+\overline{f}(\mathbf{x}_2))$ and at the same time $\overline{f}(2\mathbf{x}+2\mathbf{x}_2)=\overline{f}(2(\mathbf{x}+\mathbf{x}_2))=f(2(\mathbf{x}+\mathbf{x}_2))=2f((\mathbf{x}+\mathbf{x}_2))=2\overline{f}(\mathbf{x}+\mathbf{x}_2)$. Hence $2(\overline{f}(\mathbf{x})+\overline{f}(\mathbf{x}_2))=2\overline{f}(\mathbf{x}+\mathbf{x}_2)$ and since $T$ is torsion free we can conclude.
\end{proof}

\begin{proposition}
Let $S,T$ be GNSs in $\mathbb{N}^d$, $f:S\rightarrow T$ be an isomorphism. Then $\operatorname{g}(S)=\operatorname{g}(T)$.
\label{same-genus}
\end{proposition}
\begin{proof}
We prove by induction on $g=\operatorname{g}(S)$. If $g=0$ then $S=\mathbb{N}^d$. Suppose that $\operatorname{g}(T)\geq 1$ and consider the isomorphism $f^{-1}:T\rightarrow S$, that is the inverse of $f$. Since $\operatorname{g}(T)\geq 1$ there exists $\mathbf{x}\in \operatorname{SG}(T)$ and in particular $\mathbf{x}=\mathbf{t}-\mathbf{w}$ with $\mathbf{w}\in T\setminus \{\mathbf{0}\}$ and $\mathbf{t}\in \mathrm{Maximlas}_{\leq_T}\operatorname{Ap}(T,\mathbf{w})$. By Lemma~\ref{SG} $f^{-1}(\mathbf{t})-f^{-1}(\mathbf{w})\in \operatorname{SG}(S)$ but this is a contradiction since $S=\mathbb{N}^d$. Therefore $\operatorname{g}(T)=0$. Suppose now that $g>0$ and let $\mathbf{x}\in \operatorname{SG}(S)$. By Lemma~\ref{extend1} there exist $\mathbf{y}\in \operatorname{SG}(T)$ and an isomorphism $\overline{f}:S\cup \{\mathbf{x}\}\rightarrow T\cup \{\mathbf{y}\}$ such $\overline{f}(\mathbf{x})=\mathbf{y}$ and $\overline{f}(\mathbf{s})=f(\mathbf{s})$ for all $\mathbf{s}\in S$. Since $\operatorname{g}(S\cup \{\mathbf{x}\})=g-1$, by induction we have $\operatorname{g}(S\cup \{\mathbf{x}\})=\operatorname{g}(T\cup \{\mathbf{y}\})=\operatorname{g}(T)-1$, so $\operatorname{g}(S)=\operatorname{g}(T)$.
\end{proof}

\begin{theorem}
Let $S,T$ be GNSs in $\mathbb{N}^d$ and $f:S\rightarrow T$ be an isomorphism. Then there exists an automorphism $\varphi$ of $\mathbb{N}^d$ such that $\varphi(\mathbf{s})=f(\mathbf{s})$ for all $\mathbf{s}\in S$.
\label{extend2}
\end{theorem}
\begin{proof}
By Lemma~\ref{same-genus} we know that $\operatorname{g}(S)=\operatorname{g}(T)$. Moreover it is easy to see that there exist $\mathbf{x}_1, \mathbf{x}_2, \ldots, \mathbf{x}_n \in \operatorname{H}(S)$ such that :
\begin{itemize}
\item[a)] $\mathbf{x}_1 \in \operatorname{SG}(S)$.
\item[b)] $\mathbf{x}_i \in \operatorname{SG}(S\cup \{\mathbf{x}_1,\ldots, \mathbf{x}_{i-1}\})$ for $i\in \{2,\ldots,n\}$.
\item[c)] $S\cup \{\mathbf{x}_1,\ldots,\mathbf{x}_n\}=\mathbb{N}^d$.
\end{itemize}
By the same argument of Lemma~\ref{extend1} we can consider $\mathbf{y}_1,\mathbf{y}_2,\ldots,\mathbf{y}_n\in \operatorname{H}(T)$ such that:
\begin{itemize}
\item[a)] $\mathbf{y}_1\in \operatorname{SG}(T)$.
\item[b)] $\mathbf{y}_i\in \operatorname{SG}(T\cup \{\mathbf{y}_1,\ldots,\mathbf{y}_{i-1}\})$ for $i\in \{2,\ldots,n\}$
\item[c)] For all $i\in \{1,\ldots,n\}$, the maps $f_i:S\cup \{\mathbf{x}_1,\ldots, \mathbf{x}_{i-1}\}\rightarrow T\cup \{\mathbf{y}_1,\ldots,\mathbf{y}_{i-1}\}$ defined by $f_i(\mathbf{s})=f(\mathbf{s})$ for all $\mathbf{s}\in S$ and $f_i(\mathbf{x}_j)=\mathbf{y}_j$ for all $j\in \{1,\ldots,i\}$, are all isomorphisms.
\end{itemize}
In particular $T\cup \{\mathbf{y}_1,\ldots,\mathbf{y}_{n}\}=\mathbb{N}^d$ and $f_n:\mathbb{N}^d\rightarrow \mathbb{N}^d$ is an automorphism of $\mathbb{N}^d$ such that $f_n(\mathbf{s})=f(\mathbf{s})$ for all $\mathbf{s}\in S$.
\end{proof}

%\begin{remark}
%Let $\varphi$ be automorphism of $\mathbb{N}^d$. Consider that $\mathbb{N}^d$ is generated as a monoid by the set $\{\mathbf{e}_1,\ldots,\mathbf{e}_d\}$, that is the standard basis of the vector space $\mathbb{R}^d$. Moreover in such a case every element in $\mathbb{N}^d$ can be expressed in a unique way as $\mathbb{N}$-linear combination of the vectors $\mathbf{e}_1,\ldots,\mathbf{e}_d$. So there exists a permutation $\sigma \in \operatorname{P}_d$ such that $\varphi(\mathbf{e}_i)=\mathbf{e}_{\sigma(i)}$ for all $i\in \{1,\ldots,d\}$ and so for all $\mathbf{n}\in \mathbb{N}^d$ we have $\varphi(\mathbf{n})=\sigma(\mathbf{n})$.  
%\end{remark}

Now we can state the following result.

\begin{theorem}
Let $S,T$ be GNSs in $\mathbb{N}^d$ and suppose there exists an isomorphism $f:S\rightarrow T$. Then there exists $\sigma \in \operatorname{P}_d$ such that $T=\sigma(S)$.
\end{theorem}
\begin{proof}
Consider the isomorphism $f:S\rightarrow T$. By Theorem~\ref{extend2} then there exists an automorphism $\varphi$ of $\mathbb{N}^d$ such that $\varphi(\mathbf{s})=f(\mathbf{s})$ for all $\mathbf{s}\in S$. Consider that $\mathbb{N}^d$ is generated as a monoid by the set $\{\mathbf{e}_1,\ldots,\mathbf{e}_d\}$, that is the standard basis of the vector space $\mathbb{R}^d$. Moreover in such a case every element in $\mathbb{N}^d$ can be expressed in a unique way as $\mathbb{N}$-linear combination of the vectors $\mathbf{e}_1,\ldots,\mathbf{e}_d$. So there exists a permutation $\sigma \in \operatorname{P}_d$ such that $\varphi(\mathbf{e}_i)=\mathbf{e}_{\sigma(i)}$ for all $i\in \{1,\ldots,d\}$ and so for all $\mathbf{n}\in \mathbb{N}^d$ we have $\varphi(\mathbf{n})=\sigma(\mathbf{n})$. Hence we can easily argue that $T=f(S)=\varphi(S)=\sigma(S)$. %We prove that $T=\sigma(S)$. Let $\mathbf{t}\in T$, then there exists $\mathbf{s}=\sum_{i=1}^ds_i\mathbf{e}_i$ such that $f(\mathbf{s})=\mathbf{t}$, in particular we have $\mathbf{t}=\varphi(\mathbf{s})=\sigma(\mathbf{s})$. Moreover if $\mathbf{s}\in S$ then $\sigma(\mathbf{s})=\varphi(\mathbf{s})=f(\mathbf{s})\in T$, so we can conclude.  
\end{proof}

We saw that if $S\subseteq \mathbb{N}^d$ is a GNS then all GNSs isomorphic to $S$ are obtained considering the semigroups $\sigma(S)$ for all permutations $\sigma \in \operatorname{P}_d$. So we can introduce the equivalence relation $\simeq$ on the set $\mathcal{S}_d$ of all GNSs in $\mathbb{N}^d$, where for $S,T\in \mathcal{S}_d$ we define $S\simeq T$ if and only if $S$ and $T$ are isomorphic, equivalently there exists $\sigma\in \operatorname{P}_d$ such that $T=\sigma(S)$. For $S\in \mathcal{S}_d$, as usual we denote by $[S]_\simeq$ the equivalence class of representative $S$  with respect to the relation $\simeq$. In particular $[S]_\simeq$ is the set of all GNSs in $\mathbb{N}^d$ that are isomorphic to $S$.

\section{Equivariant generalized numerical semigroups}

We say that $S\subseteq \mathbb{N}^d$ is an \emph{equivariant} GNS if $[S]_\simeq=\{S\}$, equivalently $\sigma(S)=S$ for all $\sigma \in \operatorname{P}_d$. Another trivial equivalent condition for $S$ to be equivariant is $\sigma(\operatorname{H}(S))=\operatorname{H}(S)$ for all $\sigma \in \operatorname{P}_d$. We provide some properties of these semigroups.

\begin{example}

An example of equivariant GNS is $$S=\mathbb{N}^3\setminus \{(0, 0, 1), (0, 0, 2), (0, 1, 0), (0, 2, 0), (1, 0, 0 ), (1, 1, 1), (2, 0, 0)\}.$$
Consider also the following GNS:
$$S'=\mathbb{N}^3 \setminus \{(0, 0, 1), (0, 0, 2), (0, 0, 3), (0, 1, 0), (0, 1, 1), ( 0, 2, 0), (0, 3, 0)\}.$$
Observe that if $\sigma=(23)\in \operatorname{P}_3$, then $\sigma(S')=S'$. But $S'$ is not equivariant. In fact, if $\tau=(13)\in \operatorname{P}_3$, then $\tau(S')\neq S'$ since, for instance, $(0,0,1)\in \operatorname{H}(S')$ but $\tau((0,0,1))=(1,0,0)\notin \operatorname{H}(S')$.
%Esempio di equivariant GNS (e uno di non equivariant ma che \`e fissato da una permutazione)
\end{example}

\begin{proposition}
Let $S\subseteq \mathbb{N}^d$ be an equivariant GNS. Then the following hold:
\begin{enumerate}
\item If $\mathbf{x}$ is a minimal generator of $S$ then $\sigma(\mathbf{x})$ is a minimal generators of $S$ for all $\sigma\in \operatorname{P}_d$.
\item If $\mathbf{x}\in\operatorname{SG}(S)$ then $\sigma(\mathbf{x})\in\operatorname{SG}(S)$ for all $\sigma\in \operatorname{P}_d$.
\item If $\mathbf{x}\in\mathrm{Maximals}_{\leq}\operatorname{H}(S)$ then $\sigma(\mathbf{x})\in\mathrm{Maximals}_{\leq}\operatorname{H}(S)$ for all $\sigma\in \operatorname{P}_d$. 
\end{enumerate}
\label{3cose}
\end{proposition}
\begin{proof}
(1) Suppose there exists $\sigma \in \operatorname{P}_d$ such that $\sigma(\mathbf{x})=\mathbf{a}+\mathbf{b}$ and $\mathbf{a},\mathbf{b}\in S\setminus \{\mathbf{0}\}$. Consider $\sigma^{-1}\in \operatorname{P}_d$, then $\mathbf{x}=\sigma^{-1}(\mathbf{a})+\sigma^{-1}(\mathbf{b})$ and since $S$ is equivariant then $\sigma^{-1}(\mathbf{a}),\sigma^{-1}(\mathbf{b})\in S\setminus \{\mathbf{0}\}$. This is a contradiction, since $\mathbf{x}$ is a minimal generator of $S$.\\ 
(2) Let $\sigma \in \operatorname{P}_d$ and $\mathbf{a}\in S\setminus \{\mathbf{0}\}$. Consider $\mathbf{z}=\mathbf{a}+\sigma(\mathbf{x})$, then $\sigma^{-1}(\mathbf{z})=\sigma^{-1}(\mathbf{a})+\mathbf{x}$ and since $S$ is equivariant $\sigma^{-1}(\mathbf{a})\in S\setminus\{\mathbf{0}\}$. By hypothesis $\mathbf{x}\in\operatorname{SG}(S)$, so $\sigma^{-1}(\mathbf{z}) \in S$ and also $\mathbf{z}\in S$ since $S$ is equivariant. Moreover $2\sigma(\mathbf{x})=\sigma(2\mathbf{x})$ and considering that $\mathbf{x}\in\operatorname{SG}(S)$ and $S$ is equivariant, then $2\mathbf{x}\in S$ and so $\sigma(2\mathbf{x})=2\sigma(\mathbf{x})\in S$. Therefore $\sigma(\mathbf{x})\in \operatorname{SG}(S)$.\\
(3) Suppose there exists $\sigma \in \operatorname{P}_d$ and $\mathbf{y}\in \operatorname{H}(S)$ such that $\sigma(\mathbf{x})< \mathbf{y}$. Then it is easy to see that $\mathbf{x}< \sigma^{-1}(\mathbf{y})$ and $\sigma^{-1}(\mathbf{y})\in \operatorname{H}(S)$ since $S$ is equivariant. This is a contradiction of the maximality of $\mathbf{x}$.
\end{proof}

%\noindent Let $\mathbf{x}\in \mathbb{N}^d$, define $\mathrm{orb}(\mathbf{x})=\{\sigma(\mathbf{x})\in\mathbb{N}^d \mid \sigma \in \operatorname{P}_d\}$.

\begin{proposition}
Let $S\subseteq \mathbb{N}^d$ be an equivariant GNS. Then the following hold:
\begin{enumerate}
\item $\mathbf{x}$ is a minimal generator of $S$ if and only if $S\setminus \mathrm{orb}(\mathbf{x})$ is an equivariant GNS.
\item If $\mathbf{x}\in \mathrm{Maximals}_{\leq}\operatorname{H}(S)$ then $S\cup \mathrm{orb}(\mathbf{x})$ is an equivariant GNS.
\end{enumerate}
\label{propOrb}
\end{proposition}
\begin{proof}
(1) Suppose $\mathbf{x}$ is a minimal generator of $S$. It suffices to prove that $T=S\setminus \mathrm{orb}(\mathbf{x})$ is a semigroup. Let $\mathbf{a},\mathbf{b}\in T\setminus \{\mathbf{0}\}$, then $\mathbf{a}+\mathbf{b}\in S$. If $\mathbf{a}+\mathbf{b}=\sigma(\mathbf{x})$ for some $\sigma\in \operatorname{P}_d$ we obtain a contradiction, since $\sigma(\mathbf{x})$ is a minimal generator by Proposition~\ref{3cose}. So, $\mathbf{a}+\mathbf{b}\notin \mathrm{orb}(\mathbf{x})$, that is $\mathbf{a}+\mathbf{b}\in T$. So $T$ is a semigroup. Conversely suppose that $T=S\setminus \mathrm{orb}(\mathbf{x})$ is an equivariant GNS. If $\mathbf{x}$ is not a minimal generator of $S$ then there exists $\mathbf{a},\mathbf{b}\in S\setminus \{\mathbf{0}\}$ such that $\mathbf{x}=\mathbf{a}+\mathbf{b}$. Since $T$ and $S=T\cup \mathrm{orb}(\mathbf{x})$ are GNSs and $\mathbf{x}\notin T$ the only possibility is either $\mathbf{a}\in \mathrm{orb}(\mathbf{x})$ or $\mathbf{b}\in \mathrm{orb}(\mathbf{x})$. Without loss of generality suppose $\mathbf{a}\in \mathrm{orb}(\mathbf{x})$, in particular $\mathbf{a}=\tau(\mathbf{x})$ for some permutation $\tau \in \operatorname{P}_d$. So $\mathbf{x}=\tau(\mathbf{x})+\mathbf{b}$, that is not difficult to see it is a contradiction. Therefore $\mathbf{x}$ is a minimal generator of $S$. \\
(2) It suffices to prove that $T=S\cup \mathrm{orb}(\mathbf{x})$ is a semigroup and it easily follows considering that $\mathrm{orb}(\mathbf{x})\subseteq \mathrm{Maximals}_{\leq}\operatorname{H}(S) \subseteq \operatorname{SG}(S)$ and if $\mathbf{a},\mathbf{b}\in \mathrm{Maximals}_{\leq}\operatorname{H}(S)$ then $\mathbf{a}+\mathbf{b}\in S$. 
\end{proof}

Now we provide a construction to arrange the set of all equivariant GNSs in $\mathbb{N}^d$ in a rooted tree graph. This fact suggests a procedure to generate families of equivariant GNSs. For instance one can produce all equivariant GNSs up to a fixed value of genus.  

\begin{definition}\rm
We define $\mathcal{A}_d$ the set of all equivariant GNSs in $\mathbb{N}^d$. Let $\preceq$ be a relaxed monomial order, we define also the transform $\mathcal{F}_{\preceq}:\mathcal{A}_d\setminus \{\mathbb{N}^d\}\rightarrow \mathcal{A}_{d}$ with $\mathcal{F}_{\preceq}(S)=S\cup \mathrm{orb}(\mathbf{F}_{\preceq}(S))$.\\
Denoting by $\mathcal{F}_{\preceq}^i(S)=\mathcal{F}_{\preceq}(\mathcal{F}_{\preceq}^{i-1}(S))$, with $\mathcal{F}_{\preceq}^{1}(S)=\mathcal{F}_{\preceq}(S)$, it is not difficult to check that for all $S\in \mathcal{A}_{d}\setminus \{\mathbb{N}^d\}$ there exists $n\in \mathbb{N}$ such that $\mathcal{F}_{\preceq}^{n}(S)=\mathbb{N}^d$.

\end{definition}

\begin{lemma}
Let $S\subseteq \mathbb{N}^d$ be an equivariant GNS and $\preceq$ a relaxed monomial order. Then the following hold:
\begin{enumerate}
\item If $\mathbf{x}$ is a minimal generator of $S$ such that $\mathbf{F}_{\preceq}(S)\prec \mathbf{x}$ and $T=S\setminus \mathrm{orb}(\mathbf{x})$, then $\mathbf{F}_{\preceq}(T) \in \mathrm{orb}(\mathbf{x})$.
\item If $T=\mathcal{F}_{\preceq}(S)$ then $\mathbf{F}_{\preceq}(S)\in \mathbf{U}_{\preceq}(T)$. 
\end{enumerate}
\label{forChild}
\end{lemma}
\begin{proof}
The first claim easily follows since $\operatorname{H}(T)=\operatorname{H}(S)\cup\mathrm{orb}(\mathbf{x})$ and for all $\mathbf{h}\in \operatorname{H}(S)$ then $\mathbf{h}\prec \mathbf{x}$, so $\mathbf{F}_{\preceq}(T)=\max_{\preceq}\mathrm{orb}((\mathbf{x}))$. For the second claim, since $T\setminus \mathrm{orb}(\mathbf{F}_{\preceq}(S))=S$ is a GNS, by Proposition~\ref{propOrb} $\mathbf{F}_{\preceq}(S)$ is a minimal generator of $T$. Moreover $\operatorname{H}(T)\subset \operatorname{H}(S)$, so $\mathbf{F}_{\preceq}(T)\prec \mathbf{F}_{\preceq}(S)$. 
\end{proof}

We define the oriented graph $\mathcal{G}^{d}_{\preceq}$ having $\mathcal{A}_d$ as set of vertices and $(S,T)$ is an edge if $T=\mathcal{F}_{\preceq}(S)$.

\begin{theorem}\label{thm:equi-tree}
Let $\preceq$ be a relaxed monomial order in $\mathbb{N}^d$. Then the graph $\mathcal{G}^{d}_{\preceq}$ is a tree whose root is $\mathbb{N}^d$. Moreover let $T\in \mathcal{A}_d$ and consider the disjoint union $\mathbf{U}_{\preceq}(T)=A_{1}\sqcup A_{2}\sqcup \cdots \sqcup A_{r}$ such that if $\mathbf{x}\in A_{i}$ and $\mathbf{y}\in A_{j}$ with $i\neq j$ then $\mathbf{x}\notin \mathrm{orb}(\mathbf{y})$. For all $i\in \{1,\ldots,r\}$ consider one element $\mathbf{x}_i\in A_i$. Then the children of $T$ in $\mathcal{G}^{d}_{\preceq}$ are the semigroups $T\setminus \mathrm{orb}(\mathbf{x}_i)$, for $i \in \{1,\ldots,r\}$. 
\end{theorem}

\begin{proof}
Let $T\in \mathcal{A}_{d}$, we define the following sequence: 
\begin{itemize}
\item $T_{0}=T$.
\item $T_{i+1}=\left\lbrace\begin{array}{ll}
         \mathcal{F}_\preceq(T_{i}) & \mbox{if $T_{i}\neq \mathbb{N}^d$} \\
        \mathbb{N}^d & \mbox{otherwise}
        \end{array}
        \right.$
\end{itemize}
\noindent In particular $T_{i}=\mathcal{F}_{\preceq}^{i}(T)$ for all $i$ and there exists a non negative integer $k$ such that $T_{k}=\mathcal{F}_{\preceq}^{k}(T)=\mathbb{N}^d$. So, $(T_{0},T_{1}),(T_{1},T_{2}),\ldots,(T_{k-1},T_{k})$ provides a path from $T$ to $\mathbb{N}^d$, hence $\mathcal{G}^{d}_{\preceq}$ is a rooted tree whose root is $\mathbb{N}^d$. Let $T\in \mathcal{A}_d$ and let $S=T\setminus \mathrm{orb}(\mathbf{x}_i)$ be a semigroup such that $\mathbf{x}_i$ is an element as in the statement of this result. Since $\mathbf{F}_{\preceq}(T)\prec\mathbf{x}_i\preceq \mathbf{F}_{\preceq}(S)$ then it is easy to see that $\mathbf{F}_{\preceq}(S)\in \mathrm{orb}(\mathbf{x}_i)$, in particular $\mathrm{orb}(\mathbf{x}_i)=\mathrm{orb}(\mathbf{F}_{\preceq}(S))$. So $T=S\cup \mathrm{orb}(\mathbf{F}_{\preceq}(S))$, that is, $S$ is a child of $T$. Observe that if $\mathbf{x}_i,\mathbf{y}_i\in A_i$ then $\mathrm{orb}(\mathbf{x}_i)=\mathrm{orb}(\mathbf{y}_i)$ so it suffices to choice only one element in $A_i$ for each $i \in\{1,\ldots,r\}$. Finally if $S$ is child of $T\in \mathcal{A}_d$, by Lemma~\ref{forChild} then $T=S\cup \mathrm{orb}(\mathbf{F}_{\preceq}(S))$ with $\mathbf{F}_{\preceq}(S)\in \mathbf{U}_{\preceq}(T)$, in particular $\mathbf{F}_{\preceq}(S)\in A_i$ for some $i\in \{1,\ldots,r\}$ as in the statement of this result. Therefore we can conclude. 
\end{proof}

\begin{example}
In this example we build some branches of the tree $\mathcal{G}^{2}_{\preceq}$, where $\preceq$ is the lexicographic order. Starting from the root, $\mathbb{N}^2=\langle (0,1),(1,0)\rangle$, it is easy to see that the unique child is $R:=\mathbb{N}^2\setminus \{(0,1),(1,0)\}=\langle (0,2),(2,0),(0,3),(3,0),(1,1),(1,2),(2,1)$\}. To obtain the children of $R$, with reference to Theorem~\ref{thm:equi-tree}, we can write $\mathbf{U}_\preceq(R)=A_1 \sqcup A_2 \sqcup A_3 \sqcup A_4$, with $A_1=\{(2,0)\}$, $A_2=\{3,0\}$, $A_3=\{(1,1)\}$ and $A_4=\{(1,2),(2,1)\}$. Hence, the children of $R$ are:
\begin{itemize}
    \item $R_1=R\setminus \mathrm{orb}((2,0))=\mathbb{N}^2\setminus \{(0,1),(1,0),(0,2),(2,0)\}$.
   \item  $R_2=R\setminus \mathrm{orb}((3,0))=\mathbb{N}^2\setminus \{(0,1),(1,0),(0,3),(3,0)\}$.
   \item  $R_3=R\setminus \mathrm{orb}((1,1))=\mathbb{N}^2\setminus \{(0,1),(1,0),(1,1)\}$.
   \item  $R_2=R\setminus \mathrm{orb}((1,2))=\mathbb{N}^2\setminus \{(0,1),(1,0),(1,2),(2,1)\}$.
\end{itemize}
A picture of these branches of $\mathcal{G}^2_\preceq$ is provided in Figure~\ref{fig:ordinary-equi-tree}. Just for providing another example, we compute the children of $R_2=\mathbb{N}^2\setminus \{(0,1),(1,0),(0,3),(3,0)\}=\langle (1,1),(0,2),(2,0),(0,5),(5,0),(1,2),(2,1) \rangle$. In this case, we can write $\mathbf{U}_\preceq(R_2)=\{(5,0)\}$. So, $R_2$ has a unique child, that is, $R_2\setminus \{(0,5),(5,0)\}=\mathbb{N}^2\setminus \{(0,1),(1,0),(0,3),(3,0),(0,5),(5,0)\}$.

    \begin{figure}[h!]
\begin{tikzpicture}
\tikzset{level distance=4em}
\Tree
        [.\begin{small}$\mathbb{N}^2$\end{small}   
        \node(e){\begin{small}$\mathbb{N}^2\setminus\{(0,1),(1,0)\}$\end{small}}; 
           ]           
   
\node (e1) at (-5,-3) {\begin{small}$\mathbb{N}^2\setminus \{(0,1),(1,0),(0,2),(2,0)\}$\end{small}};
\draw (e)--(e1);
\node (e2) at (-3.2,-4.5) {\begin{small}$\mathbb{N}^2\setminus \{(0,1),(1,0),(0,3),(3,0)\}$\end{small}};
\draw (e)--(e2);
\node (e3) at (3.2,-4.5) {\begin{small}$\mathbb{N}^2\setminus \{(0,1),(1,0),(1,1)\}$\end{small}};
\draw (e)--(e3);
\node (e4) at (5,-3) {\begin{small}$\mathbb{N}^2\setminus \{(0,1),(1,0),(1,2),(2,1)\}$\end{small}};
\draw (e)--(e4);
\end{tikzpicture}
\caption{Some branches of the tree $\mathcal{G}^2_\preceq$.}
\label{fig:ordinary-equi-tree}
\end{figure}
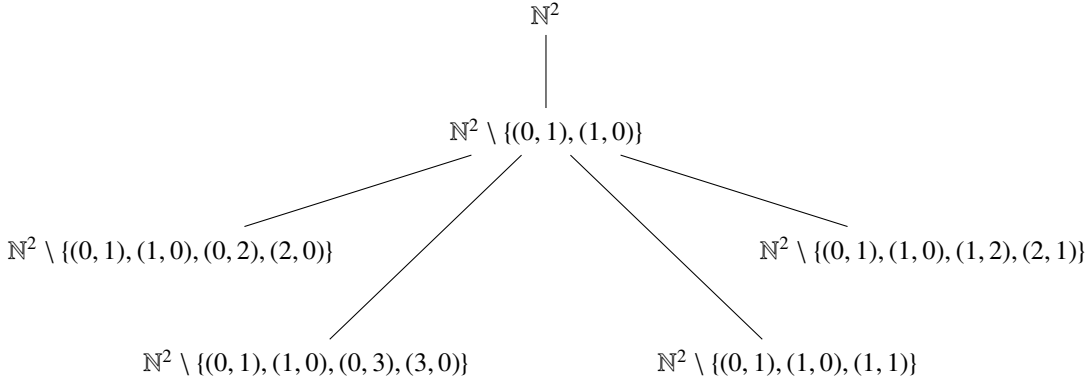

\end{example}

\section{Computing non-isomorphic generalized numerical semigroups of fixed genus}

%For $g,d$ positive integers, denote $$\mathcal{S}_{g,d}=\{S\subseteq \mathbb{N}^d\mid S\text{ is a GNS and }\operatorname{g}(S)=g\}.$$
The aim of this section is to provide some procedures to obtain for each semigroup in $\mathcal{S}_{g,d}$ only one semigroup in each equivalence class of isomorphism. That is, we want a procedure to compute the set of all non isomorphic generalized numerical semigroups in $\mathbb{N}^d$ with fixed genus $g$.

\subsection{Representatives and a first construction} The first step is to give an appropriate definition for identify a representative in $[S]_\simeq$, for any $S\in \mathcal{S}_{g,d}$. 

%\begin{definition}
%Let $\preceq$ be a total order on $\mathbb{N}^d$, we introduce an order relation $\preceq_{\operatorname{R}}$ on $\mathcal{S}_{g,d}$. Let $S,S'\in \mathcal{S}_{g,d}$, assume $\operatorname{H}(S)=\{\mathbf{h}_1 \prec \mathbf{h}_2 \prec \cdots \prec \mathbf{h}_g\}$ and $\operatorname{H}(S')=\{\mathbf{h}'_1 \prec \mathbf{h}'_2 \prec \cdots \prec \mathbf{h}'_g\}$. If $S\neq S'$, define $r=\min \{i\in \{1,\ldots,g\}\mid \mathbf{h}_i \neq \mathbf{h}'_i\}$. We set $S \preceq_{\operatorname{R}} S'$ if $S=S'$ or $\mathbf{h}_r \prec \mathbf{h}_r'$. The relation $\preceq_{\operatorname{R}}$ is trivially a total order on the set $\mathcal{S}_{g,d}$. For $S\in \mathcal{S}_{g,d}$, denote $\operatorname{R}_\preceq(S)=\min_{\leq_{\operatorname{H}}}([S]_\simeq)$ and we call it the \emph{representative} of $S$. Moreover, we define $\operatorname{R}_\preceq(\mathcal{S}_{d})=\{\operatorname{R}_\preceq(S)\mid S\in \mathcal{S}_{d}\}$ and $\operatorname{R}_\preceq(\mathcal{S}_{g,d})=\{\operatorname{R}_\preceq(S)\mid S\in \mathcal{S}_{g,d}\}$.
%\end{definition}

%$$\operatorname{r}(S,S')=1+\max \{i\in \{1,\ldots,g\}\mid \mathbf{h}_j = \mathbf{h}'_j\text{ for all }j\in \{1,\ldots,i\}\}.$$ Observe that $\operatorname{r}(S,S')=g+1$ if and only if $S=S'$. In particular, if $\operatorname{r}(S,S')\leq g$ then $\operatorname{r}(S,S')=\min \{i\in \{1,\ldots,g\}\mid \mathbf{h}_i \neq \mathbf{h}'_i\} $.

%We introduce an order relation $\preceq_{\operatorname{R}}$ on $\mathcal{S}_{g,d}$.  
\begin{definition}
Let $g,d$ be positive integers, $\preceq$ be a total order on $\mathbb{N}^d$ and $S,S'\in \mathcal{S}_{g,d}$. Assume $\operatorname{H}(S)=\{\mathbf{h}_1 \prec \mathbf{h}_2 \prec \cdots \prec \mathbf{h}_g\}$ and $\operatorname{H}(S')=\{\mathbf{h}'_1 \prec \mathbf{h}'_2 \prec \cdots \prec \mathbf{h}'_g\}$. If $S\neq S'$, define 
\[
\operatorname{r}_\prec (S,S')=\min \{i\in \{1,\ldots,g\}\mid \mathbf{h}_i \neq \mathbf{h}'_i\}
\]
 If $S=S'$ we set, by convention, $\operatorname{r}_\prec(S,S')=g+1$. We define the following relation on the set $\mathcal{S}_{g,d}$: $$S \preceq_{\operatorname{R}} S'\text{ if } S=S'\text{ or }\mathbf{h}_{\operatorname{r}_\prec(S,S')}\prec \mathbf{h}'_{\operatorname{r}_\prec(S,S')}.$$ 
The relation $\preceq_{\operatorname{R}}$ is trivially a total order on the set $\mathcal{S}_{g,d}$. For $S\in \mathcal{S}_{g,d}$, denote $\operatorname{R}_\preceq(S)=\min_{\preceq_{\operatorname{R}}}([S]_\simeq)$ and we call it the \emph{representative} of $S$ with respect to $\preceq$. Moreover, we define $\operatorname{R}_\preceq(\mathcal{S}_{d})=\{\operatorname{R}_\preceq(S)\mid S\in \mathcal{S}_{d}\}$ and $\operatorname{R}_\preceq(\mathcal{S}_{g,d})=\{\operatorname{R}_\preceq(S)\mid S\in \mathcal{S}_{g,d}\}$.
\end{definition}

\begin{example}
Consider the following GNS:
$$S=\mathbb{N}^3\setminus \{ (0,1,0),(0,2,0),(0,3,0),(1,0,0),(1,1,0),(2,0,0),(3,0,0)\}.$$
Let $\operatorname{P}_3=\{id, \sigma_2:=(12), \sigma_3:=(13), \sigma_4:=(23), \sigma_5:=(123),\sigma_6:=(132)\}$, we have:

\begin{itemize}
    \item $id(S)=\sigma_2(S)=S$.
    \item $\sigma_3(S)=\sigma_5(S)=S_2:=\mathbb{N}^3\setminus \{(0,0,1),(0,0,2),(0,0,3),(0,1,1),(0,1,0),(0,2,0),(0,3,0)\}$.
    \item $\sigma_4(S)=\sigma_6(S)=S_3:=\mathbb{N}^3\setminus \{(0,0,1),(0,0,2),(0,0,3),(1,0,0),(1,0,1),(2,0,0),(3,0,0)$\}.
\end{itemize}
Let $\preceq$ be the lexicographic order in $\mathbb{N}^3$. Then $[S]_\simeq=\{S_2 \preceq_{\operatorname{R}} S_3 \preceq_{\operatorname{R}} S\}$ and, in particular, $\operatorname{R}_\preceq(S)=S_2$.
%$$S=\mathbb{N}^3 \setminus \{(0, 0, 1), (0, 0, 2), (0, 0, 3), (0, 1, 0), (0, 1, 1), ( 0, 2, 0), (0, 3, 0)\}.$$
% \begin{itemize}
%     \item $id(S)=\sigma_4(S)=\sigma_6(S)=S$.
%     \item $\sigma_2(S)=\sigma_5(S)=S_1:=\mathbb{N}^3\setminus \{(0,0,1),(0,0,2),(0,0,3),(1,0,0),(1,0,1),(2,0,0),(3,0,0)\}$.
%     \item $\sigma_3(S)=S_2:=\mathbb{N}^3\setminus \{ (0,1,0),(0,2,0),(0,3,0),(1,0,0),(1,1,0),(2,0,0),(3,0,0)\}$. 
% \end{itemize}
\end{example}

By the previous definition, observe trivially that if $S\in \mathcal{S}_{g,d}$ is an equivariant GNS, then $\operatorname{R}_\preceq(S)=S$. Now our aim is to provide a procedure to compute directly the set $\operatorname{R}_\preceq(\mathcal{S}_{g,d})$, for $g,d$ positive integers, in such a way that the computation of non representatives is not needed. %That is, we want to provide a procedure to obtain (directly) a unique representative for each class $[S]_{\simeq}$, for every $S\in \mathcal{S}_{g,d}$. 

\begin{lemma}\label{lem:parent-is-R}
Let $S\in \operatorname{R}_\preceq(\mathcal{S}_{g,d})$ for positive integers $g,d$ and $\preceq$ be a relaxed monomial order. Then $S\cup \{\mathbf{F}_\preceq(S)\}\in \operatorname{R}_\preceq(\mathcal{S}_{g-1,d})$.
\end{lemma}
\begin{proof}
Let us denote $\mathbf{F}=\mathbf{F}_\preceq(S)$, $T=S\cup \{\mathbf{F}_\preceq(S)\}$ and assume $\operatorname{H}(S)=\{\mathbf{h}_1\prec \mathbf{h}_2\prec \cdots \prec \mathbf{h}_{g-1}\prec \mathbf{F}\}$. Trivially, we have $\operatorname{H}(T)=\{\mathbf{h}_1\prec \mathbf{h}_2\prec \cdots \prec \mathbf{h}_{g-1}\}$. Suppose $T\notin \operatorname{R}_\preceq(\mathcal{S}_{g-1,d})$. So, there exists $\sigma\in \operatorname{P}_d$ such that $\operatorname{H}(\sigma(T))=\{\sigma(\mathbf{h}_{i_1})\prec \sigma(\mathbf{h}_{i_2})\prec \cdots \prec \sigma(\mathbf{h}_{i_{g-1}})\}$ with $\{i_1,\ldots,i_{g-1}\}=\{1,\ldots,g-1\}$ and, denoted $r:=\operatorname{r}_\prec(T,\sigma(T))$, we have $\sigma(\mathbf{h}_{i_r})\prec \mathbf{h}_{r}$. % (observe that if $T=\sigma(T)$ for all $\sigma$, so $T$ is equivariant, then it belongs to $\operatorname{R}_\preceq(\mathcal{S}_{g-1,d})$). %and there exists $r=\min\{j\in \{1,\ldots,g-1\}\mid \mathbf{h}_{j}\neq \sigma(\mathbf{h}_{i_r})\}$, having $\sigma(\mathbf{h}_{i_r})\prec \mathbf{h}_{r}$ (observe that if the minimum defined by $r$ does not exists, then $T$ is equivariant and we can conclude).
Consider the element $\sigma(\mathbf{F})\in \operatorname{H}(\sigma(S))$. Observe that if $\sigma(\mathbf{h}_{i_r})\prec \sigma(\mathbf{F})$ we obtain $\sigma(S)\prec_{\operatorname{R}} S$, contradicting the assumption $S\in \operatorname{R}_\preceq(\mathcal{S}_{g,d})$. So we have $\sigma(\mathbf{F}) \prec \sigma(\mathbf{h}_{i_r})$ (observe that is not possible $\sigma(\mathbf{F})= \sigma(\mathbf{h}_{i_r})$, otherwise we obtain $\mathbf{F}=\mathbf{h}_{i_r}$). Let $k=\min\{j\in \{1,\ldots, r\}\mid \sigma(\mathbf{F})\prec \sigma(\mathbf{h}_{i_j})\}$. Observe that, since $S\prec_{\operatorname{R}} \sigma(S)$ and by our assumption on the integer $r$, we have $\mathbf{h}_k \preceq \sigma(\mathbf{F})$. If $k<r$, then $\mathbf{h}_k \preceq \sigma(\mathbf{F})\prec \sigma(\mathbf{h}_{i_k})=\mathbf{h}_k$, a contradiction. Hence $k=r$ and as a consequence, considering our assumptions, we have $\mathbf{h}_r \preceq \sigma(\mathbf{F})\prec \sigma(\mathbf{h}_{i_r})\prec \mathbf{h}_r$, obtaining again a contradiction. Therefore, we have $T\in \operatorname{R}_\preceq(\mathcal{S}_{g-1,d})$. 
\end{proof}

By Lemma~\ref{lem:parent-is-R}, we can introduce the following notions.

\begin{definition}
We introduce the transform $\mathcal{J}^{\operatorname{R}}_{\preceq}:\operatorname{R}_\preceq(\mathcal{S}_d)\setminus \{\mathbb{N}^d\}\rightarrow \operatorname{R}_\preceq(\mathcal{S}_d)$, defined by $\mathcal{J}^{\operatorname{R}}_{\preceq}(S)=S\cup \{\mathbf{F}_{\preceq}(S)\}$. Let $\mathcal{T}^{\operatorname{R}}_{d,\preceq}$  be the oriented graph having set of vertices $\operatorname{R}_\preceq(\mathcal{S}_{d})$ and $(S,T)$ is an edge of $\mathcal{T}^{\operatorname{R}}_{d,\preceq}$ if $T=\mathcal{J}^{\operatorname{R}}_{\preceq}(S)$. %In such a case we say also that $S$ is a child of $T$.
\end{definition}

\begin{theorem}\label{thm:tree-iso}
Let $\preceq$ be a relaxed monomial order in $\mathbb{N}^d$. Then the following holds:
\begin{enumerate}
\item  The graph $\mathcal{T}^{\operatorname{R}}_{d,\preceq}$ is a tree whose root is $\mathbb{N}^d$.
\item For every positive integer $g$, $S\in \operatorname{R}_\preceq(\mathcal{S}_{g,d})$ if and only if $S$ is a vertex of depth $g$ in $\mathcal{T}^{\operatorname{R}}_{d,\preceq}$.
\item For every $S\in \operatorname{R}_\preceq(\mathcal{S}_d)$, the children of $S$ in $\mathcal{T}^{\operatorname{R}}_{d,\preceq}$ are the GNSs $S\setminus \{\mathbf{n}\}$ such that $\mathbf{n}\in \mathbf{U}_\preceq(S)$ and $S\setminus \{\mathbf{n}\}\in \operatorname{R}_\preceq(\mathcal{S}_d)$.
\end{enumerate}
\end{theorem}

\begin{proof}
Let $S\in \operatorname{R}_\preceq(\mathcal{S}_{d})$, we define the following sequence: 
\begin{itemize}
\item $S_{0}=S$.
\item $S_{i+1}=\left\lbrace\begin{array}{ll}
         \mathcal{J}^{\operatorname{R}}_\preceq(S_{i}) & \mbox{if $S_{i}\neq \mathbb{N}^d$} \\
        \mathbb{N}^d & \mbox{otherwise}
        \end{array}
        \right.$
\end{itemize}
\noindent In particular $S_{i}=(\mathcal{J}_{\preceq}^{\operatorname{R}})^i(T)$ for all $i$ and it is easy to see that there exists a non negative integer $g$ such that $S_{g}=(\mathcal{J}_{\preceq}^{\operatorname{R}})^{g}(S)=\mathbb{N}^d$. So, $(S_{0},S_{1}),(S_{1},S_{2}),\ldots,(S_{g-1},S_{g})$ provides a path from $S$ to $\mathbb{N}^d$, that is, $\mathcal{T}^{\operatorname{R}}_{d,\preceq}$ is a rooted tree whose root is $\mathbb{N}^d$. Furthermore, we have that $g$ is the depth of $S$ in $\mathcal{T}^{\operatorname{R}}_{d,\preceq}$ and equivalently, by construction, we have $S=\mathbb{N}^d\setminus \{\mathbf{F}_\preceq(S_0), \mathbf{F}_\preceq(S_1), \ldots, \mathbf{F}_\preceq(S_{g-1})\}$, that is, $g=\operatorname{g}(S)$. So, we obtain the second claim. The last claim is an easy consequence of Lemma~\ref{lem:parent-is-R} and the fact that if $\mathbf{n}\in \mathbf{U}_\preceq(S)$, then $\mathbf{n}=\mathbf{F}_\preceq(S\setminus \{\mathbf{n}\})$.
\end{proof}

\begin{example}\label{exa-tree1}

We provide here an example of the construction of tree $\mathcal{T}_{d,\preceq}^{\operatorname{R}}$, where $d=2$ and $\preceq$ is the lexicographic order. Starting from $\mathbb{N}^2$, we have $\mathbf{U}_\preceq(\mathbb{N}^2)=\{(0,1),(1,0)\}$. Removing these generators, we obtain $\mathbb{N}^2\setminus \{(0,1)\}\in \operatorname{R}_\preceq(\mathcal{S}_d)$, while $\mathbb{N}^2\setminus \{(1,0)\}\notin \operatorname{R}_\preceq(\mathcal{S}_d)$. So we can consider the descendants of $\mathbb{N}^2\setminus\{(0,1)\}=\langle (0,2),(0,3),(1,0),(1,1)\rangle$. We have that $\mathbf{U}_\preceq(\mathbb{N}^2\setminus \{(0,1)\})=\{(0,2),(0,3),(1,0),(1,1)\}$, obtaining:
\begin{itemize}
    \item $S_1=\mathbb{N}^2\setminus \{(0,1),(0,2)\}\in \operatorname{R}_\preceq(\mathcal{S}_d)$.
    \item $S_2=\mathbb{N}^2\setminus \{(0,1),(0,3)\}\in \operatorname{R}_\preceq(\mathcal{S}_d)$.
    \item $S_3=\mathbb{N}^2\setminus \{(0,1),(1,0)\}\in \operatorname{R}_\preceq(\mathcal{S}_d)$.
    \item $S_4=\mathbb{N}^2\setminus \{(0,1),(1,1)\}\in \operatorname{R}_\preceq(\mathcal{S}_d)$.
\end{itemize}

\noindent Observe also that, if we consider the tree $\mathcal{T}_{d,\preceq}$ and the semigroup $S=\mathbb{N}^2\setminus \{(1,0)\}=\{(0,1),(1,1),(2,0),(3,0)\}$, having $\mathbf{U}_\preceq(S)=\{(1,1),(2,0),(3,0)\}$, we obtain:

\begin{itemize}
    \item $S_5=\mathbb{N}^2\setminus \{(1,0),(1,1)\}\notin \operatorname{R}_\preceq(\mathcal{S}_d)$ and $S_5\simeq S_4$.
    \item $S_6=\mathbb{N}^2\setminus \{(1,0),(2,0)\}\notin \operatorname{R}_\preceq(\mathcal{S}_d)$ and $S_6\simeq S_1$.
    \item $S_7=\mathbb{N}^2\setminus \{(1,0),(3,0)\}\notin \operatorname{R}_\preceq(\mathcal{S}_d)$ and $S_7\simeq S_2$.
\end{itemize}

\noindent Now, we compute the children of $S_3$ in $\mathcal{T}_{d,\preceq}^{\operatorname{R}}$. We have $S_3=\langle (0,2),(0,3),(1,1), (1,2),(2,0),(2,1),(3,0)\rangle$ and $\mathbf{U}_\preceq(S_3)=\{(1,1), (1,2), (2,0),(2,1),(3,0)\}$. So, we obtain
\begin{itemize}
    \item $T_1=\mathbb{N}^2\setminus \{(0,1),(1,0),(1,1)\}\in \operatorname{R}_\preceq(\mathcal{S}_d)$.
    \item $T_2=\mathbb{N}^2\setminus \{(0,1),(1,0),(1,2)\}\in \operatorname{R}_\preceq(\mathcal{S}_d)$.
    \item $T_3=\mathbb{N}^2\setminus \{(0,1),(1,0),(2,0)\}\notin \operatorname{R}_\preceq(\mathcal{S}_d)$.
    \item $T_4=\mathbb{N}^2\setminus \{(0,1),(1,0),(2,1)\}\notin \operatorname{R}_\preceq(\mathcal{S}_d)$.
    \item $T_5=\mathbb{N}^2\setminus \{(0,1),(1,0),(3,0)\}\notin \operatorname{R}_\preceq(\mathcal{S}_d)$.
\end{itemize}
In particular the children of $S_3$ in $\mathcal{T}_{d,\preceq}^{\operatorname{R}}$ are the semigroups $T_1$ and $T_2$. In Figure~\ref{fig:trees-s-d}, some branches of the trees $\mathcal{T}_{d,\preceq}$ and $\mathcal{T}_{d,\preceq}^{\operatorname{R}}$
are depicted, with reference to the semigroups involved in this example.
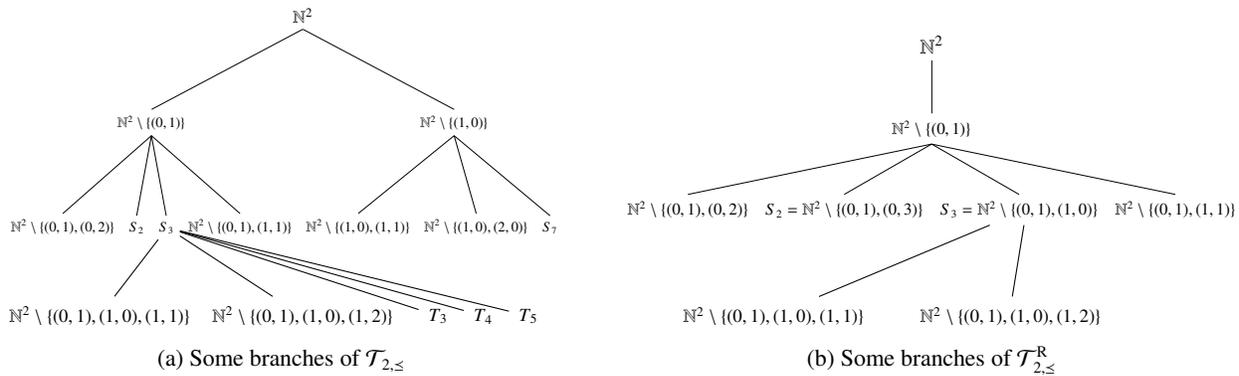
\begin{figure}[h!]
	\centering
	\subfloat[Some branches of $\mathcal{T}_{2,\preceq}$]{\begin{tikzpicture}[scale=0.6]
\tikzset{level distance=6em}
\Tree
        [.$\mathbb{N}^2$   
           [.\begin{footnotesize}$\mathbb{N}^2\setminus\{(0,1)\}$\end{footnotesize} 
                      \begin{footnotesize}$\mathbb{N}^2\setminus \{(0,1),(0,2)\}$\end{footnotesize} \begin{footnotesize}$S_2$\end{footnotesize} \node(e){\begin{footnotesize}$S_3$\end{footnotesize}};  
                      \begin{footnotesize}$\mathbb{N}^2\setminus \{(0,1),(1,1)\}$\end{footnotesize}
         ] 
          [.\begin{footnotesize}$\mathbb{N}^2\setminus \{(1,0)\}$\end{footnotesize} \begin{footnotesize}$\mathbb{N}^2\setminus \{(1,0),(1,1)\}$\end{footnotesize} \begin{footnotesize}$\mathbb{N}^2\setminus \{(1,0),(2,0)\}$\end{footnotesize} \begin{footnotesize}$S_7$\end{footnotesize} ]     
          ]
 \node (e1) at (-4.5,-6.5) {\begin{tiny}$\mathbb{N}^2\setminus \{(0,1),(1,0),(1,1)\}$\end{tiny}};
\draw (e)--(e1);
\node (e2) at (3,-6.5) {\begin{tiny}$T_3$\end{tiny}};
\draw (e)--(e2);
\node (e3) at (4,-6.5) {\begin{tiny}$T_4$\end{tiny}};
\draw (e)--(e3);
\node (e4) at (0,-6.5) {\begin{tiny}$\mathbb{N}^2\setminus \{(0,1),(1,0),(1,2)\}$\end{tiny}};
\draw (e)--(e4);
\node (e5) at (5,-6.5) {\begin{tiny}$T_5$\end{tiny}};
\draw (e)--(e5);
\end{tikzpicture}}
	\qquad
	\subfloat[Some branches of $\mathcal{T}_{2,\preceq}^{\operatorname{R}}$]{\begin{tikzpicture}[scale=0.7]
\tikzset{level distance=4em}
\Tree
        [.$\mathbb{N}^2$   
           [.\begin{footnotesize}$\mathbb{N}^2\setminus\{(0,1)\}$\end{footnotesize} 
                      \begin{footnotesize}$\mathbb{N}^2\setminus \{(0,1),(0,2)\}$\end{footnotesize} \begin{footnotesize}$S_2=\mathbb{N}^2\setminus \{(0,1),(0,3)\}$\end{footnotesize} \node(e){\begin{footnotesize}$S_3=\mathbb{N}^2\setminus \{(0,1),(1,0)\}$\end{footnotesize}};  \begin{footnotesize}$\mathbb{N}^2\setminus \{(0,1),(1,1)\}$\end{footnotesize}
         ]              
   ]
   
\node (e1) at (-3,-5) {\begin{tiny}$\mathbb{N}^2\setminus \{(0,1),(1,0),(1,1)\}$\end{tiny}};
\draw (e)--(e1);
\node (e4) at (1.5,-5) {\begin{tiny}$\mathbb{N}^2\setminus \{(0,1),(1,0),(1,2)\}$\end{tiny}};
\draw (e)--(e4);
\end{tikzpicture}} 
	\caption{Some branches of $\mathcal{T}_{2,\preceq}$ and $\mathcal{T}_{d,\preceq}^{\operatorname{R}}$ where $\preceq$ is the lexicographic order, with reference to Example~\ref{exa-tree1}.}
	\label{fig:trees-s-d}
\end{figure}

\end{example}

\begin{remark}\label{rmk:trees}
Let $g,d$ be a positive integers and $\preceq$ a relaxed monomial in $\mathbb{N}^d$. Observe that, roughly speacking, Theorem~\ref{thm:tree-iso} is nothing more than the restriction on $\operatorname{R}_{\preceq}(\mathcal{S}_d)$ of the construction on $\mathcal{S}_d$ of the semigroup tree $\mathcal{T}_{d,\preceq}$. This construction is allowed, in a certain sense, since the transform $S\mapsto S\cup \{\mathbf{F}_\preceq(S)\}$ can be restricted from $\mathcal{S}_d$ to $\operatorname{R}_\preceq(\mathcal{S}_d)$, by Lemma~\ref{lem:parent-is-R}.
\end{remark}

\begin{algorithm} 
\DontPrintSemicolon

\KwData{Two integers $g,d\in \mathbb{N}$ and a relaxed monomial order $\preceq$.}
\KwResult{$\operatorname{R}_\preceq(\mathcal{S}_{g,d})$}

 $G=\{\mathbf e_1, \mathbf e_2, \dots, \mathbf e_d\}$, $R_{0,d}=\{\mathbb{N}^{d}\}$, $\mathbf{F}_{\preceq}(\mathbb{N}^{d})=(-1,\ldots,-1)$, $\operatorname{H}(\mathbb{N}^{d})=\emptyset$.\; \label{1}
    \For{$i\in\{0,\ldots,g\}$}{
 Denote $R_{i,d}=\{S^{(j)}\mid j \in \{1,\ldots,|R_{i,d}|\}\}$.\;
$R_{i+1,d}=\emptyset$.\;
     \For{$j\in \{1,\ldots, |R_{i,d}|\}$}{
     From $S^{(j)}$ compute $E^{(j)}=\mathbf{U}_\preceq(S^{(j)})$.\;
     Denote $E^{(j)}=\{\mathbf{g}_{1},\mathbf{g}_{2},\ldots,\mathbf{g}_{|E^{(j)}|}\}$.\;
            \For{$k\in\{1,\ldots,|E^{(j)}|\}$}{
              $S^{(j,k)}=S^{(j)}\setminus \{\mathbf{g}_{k}\}$\;            
            
             \nl \If{$S^{(j,k)}=\min_{\preceq_{\operatorname{R}}}([S^{(j,k)}]_\simeq)$ \label{1}}{ 
                $S_{i+1,d}=S_{i+1,d}\cup \{(S^{(j,k)}\}$.\;
                $\mathbf{F}_{\preceq}(S^{(j,k)})=\mathbf{g}_{k}$.\;
              \nl  $\operatorname{H}(S^{(j,k)})=\operatorname{H}(S^{(j)})\cup \{\mathbf{g}_{k}\}$.\; \label{2}    
               } } 
               \If{$i+1=g$}{\Return $R_{g,d}$.\;}      
            }       
     }
\caption{Algorithm for computing $\operatorname{R}_\preceq(\mathcal{S}_{g,d})$} \label{alg:ngd}
\end{algorithm}

\noindent A procedure to compute the set $\operatorname{R}_\preceq(\mathcal{S}_{g,d})$ using the construction suggested by Theorem~\ref{thm:tree-iso} is provided in Algorithm~\ref{alg:ngd}. The instruction in line \ref{2} is useful to update the set of gaps of each semigroup computed. This fact allows to avoid the computation of the set of gaps from scratch (the known algorithms use the minimal generators for this) in the case this set could be used inside the algorithm (for instance, as we explain below, for the condition in line \ref{1}) or one is interested to have the sets of gaps at end of the algorithm. 
The instruction in line \ref{1} is crucial: only the semigroups that are representatives in their equivalence class must be used to continue the procedure. The direct way is to compute the action of all permutations in the set of gaps of $S^{(j)}$ and then compare them with respect to $\preceq_{\operatorname{R}}$.  For a more efficient implementation, it would be useful an efficient computationally procedure to check if $S^{(j,k)}=\min_{\preceq_{\operatorname{R}}}([S^{(j,k)}]_\simeq)$. More in general, if $T\in \operatorname{R}_\preceq(\mathcal{S}_{g-1,d})$ and $\mathbf{n}$ is a minimal generator of $T$, it would be useful to find conditions in order to $T\setminus \{\mathbf{n}\}\in \operatorname{R}_\preceq(\mathcal{S}_{g,d})$. The next results allows us to provide a sufficient condition for this.

%\begin{lemma}\label{lem:R-remove-gen1}
%Let $g,d$ be positive integers, $\preceq$ be a relaxed monomial order and $T\in \operatorname{R}_\preceq(\mathcal{S}_{g-1,d})$. Let $\mathbf{n}$ be a minimal generator of $T$ and $\sigma\in\operatorname{P}_d$ such that $\mathbf{n}\preceq \sigma(\mathbf{n})$. Then $T\setminus \{\mathbf{n}\}\preceq_{\operatorname{R}} \sigma(T\setminus \{\mathbf{n}\})$.
%\end{lemma}

\begin{lemma} \label{lem:minimum}
    Let $d$ be positive integer and $\preceq$ be a relaxed monomial order. If $S\in \operatorname{R}_\preceq(\mathcal{S}_d)$ then 
    $$\min_\preceq \operatorname{H}(S)=\min_\preceq \mathbb{N}^d\setminus \{\mathbf{0}\}=\min_\preceq (\{\mathbf{e}_1,\ldots,\mathbf{e}_d\}).$$
\end{lemma}
\begin{proof}
    Observe that for every $\mathbf{n}\in \mathbb{N}^d$ there exists $i\in \{1,\ldots,d\}$ such that $\mathbf{e}_i\leq \mathbf{n}$. Since every relaxed monomial order extends the natural partial in $\mathbb{N}^d$ (see \cite[Proposition 4.4]{cisto2019irreducible}), we have $\mathbf{e}_i\preceq \mathbf{n}$. So, the second equality in the statement holds. We prove that $\min_\preceq \operatorname{H}(S)=\min_\preceq (\{\mathbf{e}_1,\ldots,\mathbf{e}_d\}).$ Assume $\min_\preceq \operatorname{H}(S)=\mathbf{h}$. We first show that $\mathbf{h}\in \{\mathbf{e}_1,\ldots,\mathbf{e}_d\}$. It is sufficient to prove that there exists $\mathbf{e}_\ell$ such that $\mathbf{e}_\ell \leq \mathbf{h}$ and $\mathbf{e}_\ell\in \operatorname{H}(S)$, since in this case we have $\mathbf{e}_\ell\preceq \mathbf{h}$ and $\mathbf{h}\preceq \mathbf{e}_\ell$ by the minimality of $\mathbf{h}$. It is not difficult to see that there exists $\mathbf{e}_i$ such that $\mathbf{h}=\mathbf{h}'+\mathbf{e}_i$ for some $i\in \{1,\ldots,d\}$ and $\mathbf{h}'\in \mathbb{N}^d$. Since $\mathbf{h}\in \operatorname{H}(S)$, we have $\mathbf{e}_i\in \operatorname{H}(S)$ or $\mathbf{h}'\in \operatorname{H}(S)$. If $\mathbf{e}_i\in\ \operatorname{H}(S)$ we obtain our claim. If $\mathbf{e}_i\notin \operatorname{H}(S)$, then $\mathbf{h}'\in \operatorname{H}(S)$ and there exists $\mathbf{e}_j$ such that $\mathbf{e}_j\leq \mathbf{h}'$. In particular, $\mathbf{e}_j\in \operatorname{H}(S)$ or $\mathbf{h}'-\mathbf{e}_j\in \operatorname{H}(S)$.
    Continuing this procedure, at the end we obtain that there exists $\mathbf{e}_\ell$ such that $\mathbf{e}_\ell \leq \mathbf{h}$ and $\mathbf{e}_\ell\in \operatorname{H}(S)$. Therefore, $\mathbf{h}=\mathbf{e}_\ell$ for some $\ell\in \{1,\ldots,d\}$. Let $\min_\preceq(\{\mathbf{e}_1,\ldots,\mathbf{e}_d\})=\mathbf{e}_k$ and suppose $\mathbf{e}_\ell\neq \mathbf{e}_k$. Let $\sigma\in \operatorname{P}_d$ such that $\sigma(\mathbf{e}_\ell)=\mathbf{e}_k$. Then we obtain $\min_\preceq \operatorname{H}(\sigma(S))=\mathbf{e}_k$ and $\min_\preceq \operatorname{H}(S)=\mathbf{e}_\ell$ with $\mathbf{e}_k \prec \mathbf{e}_\ell$, that contradicts $S\in \operatorname{R}_\preceq(\mathcal{S}_d)$. 
\end{proof}

% \begin{proposition}\label{lem:R-remove-gen1}
% Let $g,d$ be positive integers, $\preceq$ be a relaxed monomial order and $T\subseteq \mathbb{N}^d$ be a GNS. Let $\sigma\in\operatorname{P}_d$ such that $T\preceq_{\operatorname{R}} \sigma(T)$. Assume $\mathbf{n}$ is a minimal generator of $T$ such that $\mathbf{n}\preceq \sigma(\mathbf{n})$. Then $T\setminus \{\mathbf{n}\}\preceq_{\operatorname{R}} \sigma(T\setminus \{\mathbf{n}\})$.
% \end{proposition}

\begin{proposition}\label{prop:R-remove-gen1}
Let $g,d$ be positive integers, $\preceq$ be a relaxed monomial order and $T\in \operatorname{R}_\preceq(\mathcal{S}_{g-1,d})$. Let $\mathbf{n}$ be a minimal generator of $T$ such that $\mathbf{n}\preceq \sigma(\mathbf{n})$ for all $\sigma \in \operatorname{P}_d$. Then $T\setminus \{\mathbf{n}\}\in \operatorname{R}_\preceq(\mathcal{S}_{g,d})$.
\end{proposition}

\begin{proof}
Let $\sigma \in \operatorname{P}_d$, we want to prove that $T\setminus \{\mathbf{n}\}\preceq_{\operatorname{R}} \sigma(T\setminus \{\mathbf{n}\})$. We know that $T\setminus \{\mathbf{n}\}$ is a GNS. Assume $T$ has genus $g-1$ for some integer $g>1$. Denote $\operatorname{H}(T)=\{\mathbf{h}_1 \prec \mathbf{h}_2 \prec \cdots \prec \mathbf{h}_{g-1} \}$ and $\operatorname{H}(\sigma(T))=\{\sigma(\mathbf{h}_{i_1})\prec \sigma(\mathbf{h}_{i_2})\prec \cdots \prec \sigma(\mathbf{h}_{i_{g-1}})\}$ with $\{i_1,\ldots,i_{g-1}\}=\{1,\ldots,g-1\}$. Denote also $r:=\operatorname{r}_\prec(T,\sigma(T))$, in particular we have either $T=\sigma(T)$ or $\mathbf{h}_{r}\prec \sigma(\mathbf{h}_{i_r})$, and $\mathbf{h}_j=\sigma(\mathbf{h}_{i_j})$ for all $j\in \{1,\ldots,r-1\}$ in the case $r>1$. 
%Observe that if $\mathbf{h}_{g-1}\prec \mathbf{n}$, then we can easily conclude. So, we can assume there exists $k=\min\{j\in \{1,\ldots,g-1\}\mid \mathbf{n}\prec \mathbf{h}_j\}$. If $\sigma(\mathbf{n})\prec \sigma(\mathbf{h}_{i_{g-1}})$, we define the integer $\ell=\min \{j\in \{1,\ldots,g\}\mid \sigma(\mathbf{n})\prec \sigma(\mathbf{h}_{i_j})\}$, while we set $\ell=g$ if $\sigma(\mathbf{h}_{i_{g-1}})\prec \sigma(\mathbf{n})$. 
Observe that, by Lemma~\ref{lem:minimum}, $\mathbf{h}_1\prec \mathbf{n}$, since $\mathbf{n}\in T$. Hence, we can define the integer $k=1+\max\{j\in \{1,\ldots,g-1\}\mid \mathbf{h}_j \prec \mathbf{n}\}$.
Moreover, we can assume that the integer $\ell=1+\max \{j\in \{1,\ldots,g-1\}\mid  \sigma(\mathbf{h}_{i_j})\prec \sigma(\mathbf{n}) \}$ is well defined. Otherwise, we have $\sigma(\mathbf{n})\prec \sigma(\mathbf{h}_{i_1})$ and, by Lemma~\ref{lem:minimum}, we have $\mathbf{h}_1\preceq \sigma(\mathbf{n})$. If $\mathbf{h}_1\prec \sigma(\mathbf{n})$ we trivially obtain $T\setminus \{\mathbf{n}\}\preceq_{\operatorname{R}} \sigma(T\setminus \{\mathbf{n}\})$. If $\mathbf{h}_1= \sigma(\mathbf{n})$, since $\mathbf{n}\preceq \sigma(\mathbf{n})$, using again Lemma~\ref{lem:minimum} the only possibility is $\mathbf{n}=\mathbf{h}_1$, a contradiction. About the integers $k,\ell\in \{2,\ldots,g\}$, we distinguish the following cases:
\begin{enumerate}
\item $k=\ell$. In this case, it is not difficult to see that   $T\setminus \{\mathbf{n}\}\preceq_{\operatorname{R}} \sigma(T\setminus \{\mathbf{n}\})$, and we can conclude.
\item $\ell < k$ and $\ell < r$. In this case we have $\sigma(\mathbf{n})\prec \sigma(\mathbf{h}_{i_\ell})=\mathbf{h}_{\ell}\preceq \mathbf{h}_{k-1}\prec \mathbf{n}$, that contradicts $\mathbf{n}\preceq \sigma(\mathbf{n})$. So, this case cannot occur with our hypothesis.
\item $r<k$ and $r<\ell$. In this case, since $\mathbf{h}_{r}\prec \sigma(\mathbf{h}_{i_r})$, with $\mathbf{h}_j=\sigma(\mathbf{h}_{i_j})$ for all $j\in \{1,\ldots,r-1\}$ in the occurrence $r>1$, it is not difficult to see that $T\setminus \{\mathbf{n}\}\prec_{\operatorname{R}} \sigma(T\setminus \{\mathbf{n}\})$.
\item $r=\ell<k$. In this case observe that, by construction, $\mathbf{h}_r \prec \mathbf{n}$. If we suppose that $\sigma(\mathbf{n})\preceq \mathbf{h}_{r}$, we obtain $\sigma(\mathbf{n})\prec \mathbf{n}$. But this contradicts our hypothesis. So, $\mathbf{h}_r \prec \sigma(\mathbf{n})$. In this case, this means $T\setminus \{\mathbf{n}\}\prec_{\operatorname{R}} \sigma(T\setminus \{\mathbf{n}\})$. %Assuming%$T\setminus \{\mathbf{n}\}\prec_{\operatorname{R}} \sigma(T\setminus \{\mathbf{n}\})$
%\item $r=k<\ell$. In this case, by construction, we have $\mathbf{n}\prec \mathbf{h}_r \prec \sigma(\mathbf{h_{i_r}})$. In particular, we have $\mathbf{n}\prec \sigma(\mathbf{h_{i_r}})$, that means $T\setminus \{\mathbf{n}\}\prec_{\operatorname{R}} \sigma(T\setminus \{\mathbf{n}\})$.
%\item $k\leq r<\ell$.
\item $k\leq r$ and $k<\ell$. In this case, by construction, we have $\mathbf{n}\prec \mathbf{h}_k \preceq \sigma(\mathbf{h}_{i_k})$. In particular, we have $\mathbf{n}\prec \sigma(\mathbf{h}_{i_k})$, $\mathbf{h}_{k-1}\prec \mathbf{n}$ and $\mathbf{h}_j=\sigma(\mathbf{h}_{i_j})$ for all $j\in \{1,\ldots,k-1\}$. Hence, $T\setminus \{\mathbf{n}\}\prec_{\operatorname{R}} \sigma(T\setminus \{\mathbf{n}\})$.
%\item $k<\ell \leq r$. In this case, $\mathbf{n}\prec \mathbf{h}_k=\sigma(\mathbf{h}_{i_k})$ and if $k>1$ we have also $\mathbf{h}_{k-1}\prec \mathbf{n}$ and $\mathbf{h}_j=\sigma(\mathbf{h}_{i_j})$ for all $j\in \{1,\ldots,k-1\}$. Hence, $T\setminus \{\mathbf{n}\}\prec_{\operatorname{R}} \sigma(T\setminus \{\mathbf{n}\})$.
\end{enumerate}
The previous points describe all possible occurrences. So, we can conclude.
\end{proof}

%Observe that if $\mathbf{h}_{g-1}\prec \mathbf{n}$, then we can easily conclude. So, we can assume there exists $k=\min\{j\in \{1,\ldots,g-1\}\mid \mathbf{n}\prec \mathbf{h}_j\}$. If $\sigma(\mathbf{n})\prec \sigma(\mathbf{h}_{i_{g-1}})$, we define the integer $\ell=\min \{j\in \{1,\ldots,g\}\mid \sigma(\mathbf{n})\prec \sigma(\mathbf{h}_{i_j})\}$, while we set $\ell=g$ if $\sigma(\mathbf{h}_{i_{g-1}})\prec \sigma(\mathbf{n})$. ALTERNATIVA A SOTTO

% \begin{corollary}\label{cor:R-remove-gen1}
% Let $g,d$ be positive integers, $\preceq$ be a relaxed monomial order and $T\in \operatorname{R}_\preceq(\mathcal{S}_{g-1,d})$. Let $\mathbf{n}$ be a minimal generator of $T$ such that $\mathbf{n}\preceq \sigma(\mathbf{n})$ for all $\sigma \in \operatorname{P}_d$. Then $T\setminus \{\mathbf{n}\}\in \operatorname{R}_\preceq(\mathcal{S}_{g,d})$.
% \end{corollary}

%\textcolor{red}{\textbf{Question}: Is there a characterization for a minimal generator $\mathbf{n}\in T$ in order to $T\setminus \{\mathbf{n}\}\in \operatorname{R}_\preceq(\mathcal{S}_{g,d})$? Is there a characterization at least for effective generators?}

%\subsection{Some refinements for the previous construction}\textcolor{red}{TO BE WRITTEN}

 The next result states another very particular condition for a semigroup $T\in \operatorname{R}_\preceq(\mathcal{S}_d)$ and $\mathbf{n}\in \mathbf{U}_\preceq(T)$ such that $T\setminus \{\mathbf{n}\}\in \operatorname{R}_\preceq(\mathcal{S}_d)$, in the case $d=2$. If $S\in \operatorname{R}_\preceq(\mathcal{S}_d)$ we denote by $\mathcal{D}(S)$ the set containing $S$ and all descendants of $S$ in the tree $\mathcal{T}_{d,\preceq}^{\operatorname{R}}$. In particular, if $T\in \mathcal{D}(S)$, then there exists a sequence of GNSs $T_1,\ldots,T_n\in \operatorname{R}_\preceq(\mathcal{S}_d)$, such that $T_1=T$, $T_n=S$ and $T_{i+1}=T_{i}\cup \{\mathbf{F}_\preceq (T_i)\}$ for all $i\in \{1,\ldots,n-1\}$. 

\begin{proposition}\label{prop:descendant_N2}
    Let $\preceq$ be a relaxed monomial order in $\mathbb{N}^2$ and assume $\mathbf{e}_1\prec \mathbf{e}_2$. Let $S\in \operatorname{R}_\preceq(\mathcal{S}_2)\setminus \{\mathbb{N}^2\}$ such that $\mathbf{e}_2 \in S$ and $\mathbf{e}_2\prec \mathbf{F}_\preceq (S)$.
    Then for all $T\in \mathcal{D}(S)$ and for all $\mathbf{n}\in \mathbf{U}_\preceq(T)$ we have $T\setminus \{\mathbf{n}\}\in \operatorname{R}_\preceq(\mathcal{S}_2)$
\end{proposition}
\begin{proof}
     By hypothesis, if $T\in \mathcal{D}(S)$, then it is not difficult to see that $\mathbf{e}_1\in \operatorname{H}(T)$ (also by Lemma~\ref{lem:minimum}), $\mathbf{e}_2 \notin \operatorname{H}(T)$ and $\mathbf{e}_2\prec \mathbf{F}_\preceq(T)$. Let $\mathbf{n}\in \mathbf{U}_\preceq (T)$, we show that $T'=T\setminus \{\mathbf{n}\}\in \operatorname{R}_\preceq(\mathcal{S}_2)$. Observe that also for $T'$ we have $\mathbf{e}_1 \operatorname{H}(T')$, $\mathbf{e}_2 \notin \operatorname{H}(T')$ and $\mathbf{e}_2\prec \mathbf{F}_\preceq(T')$. In particular, if $\sigma$ is the permutation in $\operatorname{P}_2$ such that $\sigma\neq \mathrm{id}$, then $\mathbf{e}_1\notin \sigma(T')$ and by Lemma~\ref{lem:minimum} we have $\sigma(T')\notin \operatorname{R}_\preceq(\mathcal{S}_2)$. 
\end{proof}

If $\preceq$ is the lexicographic order, an example of GNS satisfying Proposition~\ref{prop:descendant_N2} is $S=\mathbb{N}^2 \setminus \{(0,1),(1,1)\}$. Moreover, observe that if in Proposition~\ref{prop:descendant_N2} we assume that $\preceq$ is the graded lexicographic order, $\mathbf{e}_2\in S$ and $\operatorname{g}(S)\leq 2$, then the assumption $\mathbf{e}_2\prec \mathbf{F}_\preceq (S)$ is automatically verified. Inspired by some example we computed, in the following we provide other properties of the semigroups in $\operatorname{R}_\preceq(\mathcal{S}_d)$, which generalize Lemma~\ref{lem:minimum}.

\begin{remark}\label{rem:trasposizioni}
    Let $d$ be a positive integer and $\preceq$ be a relaxed monomial order. Suppose that $\mathbf{e}_1\prec \mathbf{e}_2\prec \cdots \prec \mathbf{e}_d$ in $\mathbb{N}^d$. Let $\{\mathbf{e}_{j_1}\prec \mathbf{e}_{j_2}\prec\ldots\prec\mathbf{e}_{j_r}\}\subseteq \{\mathbf{e}_1,\ldots,\mathbf{e}_d\}$ for some $\{j_1,j_2,\ldots,j_r\}\subseteq \{1,\ldots,d\}$ and $r\leq d$. Then there exists a permutation $\sigma\in \operatorname{P}_d$ such that $\sigma(\{\mathbf{e}_{j_1}, \mathbf{e}_{j_2},\ldots,\mathbf{e}_{j_r}\})=\{\mathbf{e}_1,\ldots,\mathbf{e}_r\}$. In fact, it is sufficient to take $\sigma=(j_rr)\ldots(j_22)(j_11)$, as product of disjoint transpositions (and conventionally, $(j_k k)=\mathrm{id}$ if $j_k=k$).
\end{remark}

%   \noindent  Moreover, with the same assumptions, if we consider the set $E=\{\mathbf{e}_1,\ldots,\mathbf{e}_{d-1}\}$, then for all $\sigma\in \operatorname{P}_d\setminus \{\mathrm{id}\}$ we have $\sigma(E)\neq E$. \textcolor{red}{to add a brief explanation?} QUESTA COMMENTATA ERA SBAGLIATISSIMA, BASTA PRENDERE LA PERMUTAZIONE (2 3 ... n)

\begin{proposition}\label{prop:basis-2}
    Let $d$ be a positive integer and $\preceq$ be a relaxed monomial order. Suppose that $\mathbf{e}_1\prec \mathbf{e}_2\prec \cdots \prec \mathbf{e}_d$ in $\mathbb{N}^d$. Let $S\in \operatorname{R}_\preceq(\mathcal{S}_d)$ and suppose $\operatorname{H}(S)\cap \{\mathbf{e}_1,\ldots,\mathbf{e}_d\}=\{\mathbf{e}_{j_1}\prec \mathbf{e}_{j_2}\}$ for some $\{j_1,j_2\}\subseteq \{1,\ldots,d\}$. Then $\{\mathbf{e}_{j_1}, \mathbf{e}_{j_2}\}=\{\mathbf{e}_{1},\mathbf{e}_{2}\}$. 
\end{proposition}
\begin{proof}
    By Lemma~\ref{lem:minimum} we have $j_1=1$ and we can assume $\operatorname{H}(S)=\{\mathbf{e}_1 \prec \mathbf{h}_1 \prec \cdots \prec \mathbf{h}_t \prec \mathbf{e}_{j_2}\prec \cdots \}$. Let $T=\mathbb{N}^d \setminus \{\mathbf{e}_1 \prec \mathbf{h}_1 \prec \cdots \prec \mathbf{h}_t \prec \mathbf{e}_{j_2}\}$. By Lemma~\ref{lem:parent-is-R} we have $T\in \operatorname{R}_\preceq (\mathcal{S}_d)$. Let $\sigma \in \operatorname{P}_d$ such that $\sigma(\{\mathbf{e}_1,\mathbf{e}_{j_2}\})=\{\mathbf{e}_1,\mathbf{e}_2\}$ and consider $\operatorname{H}(\sigma(T))=\{\mathbf{e}_1\prec \mathbf{g}_1\prec \cdots \prec \mathbf{g}_m \prec \mathbf{e}_2 \prec \mathbf{g}_{m+1}\prec \cdots \mathbf{g}_t\}$ with $m\leq t$ (where $m=t$ if and only if $\operatorname{H}(\sigma(T))=\{\mathbf{e}_1\prec \mathbf{g}_1\prec \cdots \prec \mathbf{g}_t \prec \mathbf{e}_2\}$). Observe that, for $i\leq m$, we have $\mathbf{g}_i\prec \mathbf{e}_j$ for all $j\geq 2$. So, since every relaxed monomial order extends the natural partial order in $\mathbb{N}^d$, for $i\leq m$ the only possibility is $\mathbf{g}_i=\alpha_i\mathbf{e}_1$ for some $\alpha_i\in\mathbb{N}\setminus \{0,1\}$. Define $\mathbf{h}_{k_i}$ such that $\sigma(\mathbf{h}_{k_i})=\mathbf{g}_i$. In particular, $\mathbf{h}_{k_i}=\sigma^{-1}(\mathbf{g}_i)=\alpha_i \sigma^{-1}(\mathbf{e}_1)$. We can consider two possibility: $\sigma^{-1}(\mathbf{e}_1)=\mathbf{e}_{j_2}$ or $\sigma^{-1}(\mathbf{e}_1)=\mathbf{e}_1$. In the first case, we obtain $\mathbf{h}_{k_i}=\alpha_i \mathbf{e}_{j_2}\succ \mathbf{e}_{j_2}$, that is a contradiction. So, we have $\sigma^{-1}(\mathbf{e}_1)=\mathbf{e}_1$, that implies $\sigma(\mathbf{e}_1)=\mathbf{e}_1$. This means that $\sigma(\mathbf{h}_{k_i})=\mathbf{h}_{k_i}=\mathbf{g}_i=\alpha_i\mathbf{e}_1$ for all $i\leq m$. In particular, we have also $\mathbf{h}_{k_1}\prec \mathbf{h}_{k_2}\prec \cdots \prec \mathbf{h}_{k_m}$. We want to show that $\mathbf{h}_{k_i}=\mathbf{h}_i$ for all $i\leq m$. In fact, assume for instance that $\mathbf{h}_1\notin \{\mathbf{h}_{k_1},\ldots,\mathbf{h}_{k_m}\}$. As a consequence, we have $\sigma(\mathbf{h}_1)=\mathbf{g}_k$ for some $k>m$. In particular $\mathbf{e}_2\prec \sigma(\mathbf{h}_1)$. Another consequence is $\mathbf{e}_1\prec \mathbf{h}_1\prec \mathbf{h}_{k_1}\prec \mathbf{e}_2$, from which we obtain $\mathbf{h}_1=\beta \mathbf{e}_1$, for some $\beta\in\mathbb{N}\setminus \{0,1\}$. So, $\sigma(\mathbf{h}_1)=\beta \sigma(\mathbf{e}_1)=\beta \mathbf{e}_1=\mathbf{h}_1\prec \mathbf{e}_2$, a contradiction. Hence, $\mathbf{h}_1\in \{\mathbf{h}_{k_1},\ldots,\mathbf{h}_{k_m}\}$ and the only possibility is $\mathbf{h}_1=\mathbf{h}_{k_1}$.  By the same argument, if we assume $\mathbf{h}_2\notin \{\mathbf{h}_{k_2},\ldots,\mathbf{h}_{k_m}\}$, then we obtain obtain $\mathbf{h}_2=\mathbf{h}_{k_2}$, and continuing we have $\mathbf{h}_i=\mathbf{h}_{k_i}$ for all $i\leq m$. Therefore, $\mathbf{h}_i=\sigma(\mathbf{h}_i)=\mathbf{g}_i$ for all $i\leq m$. Now, suppose $m<t$. Considering the element $\mathbf{h}_{m+1}$, by the previous arguments the only possibility is $\sigma(\mathbf{h}_{m+1})=\mathbf{g}_k$ for $k> m$, that is $\mathbf{e}_2\prec \sigma(\mathbf{h}_{m+1})$. Moreover, by $T\in \operatorname{R}_\preceq(\mathcal{S}_d)$, we have $\mathbf{h}_{m+1}\prec \mathbf{e}_2$. In particular, $\mathbf{h}_{m+1}=\gamma\mathbf{e}_1$ for some $\gamma\in \mathbb{N}\setminus \{0,1\}$, that implies $\sigma(\mathbf{h}_{m+1})=\mathbf{h}_{m+1}\prec \mathbf{e}_2$, a contradiction. So, $m=t$ and in order to $T\in \operatorname{R}_\preceq(\mathcal{S}_d)$ the only possibility is $\mathbf{e}_{j_2}\preceq \mathbf{e}_2$. By construction, this implies $\mathbf{e}_{j_2}=\mathbf{e}_2$.
\end{proof}

The next property can be seen as  a generalization of Proposition~\ref{prop:basis-2}, considering more than two standard basis vectors, but in the case the fixed relaxed monomial order has a particular behaviour.

\noindent If $\preceq$ is a relaxed monomial order in $\mathbb{N}^d$, we say that $\preceq$ is \emph{1-graded} if for all $\mathbf{n}\in \mathbb{N}^d\setminus \{\mathbf{0}, \mathbf{e}_1,\ldots,\mathbf{e}_d\}$ it is verified $\mathbf{n}\succ \max_\preceq \{\mathbf{e}_1,\ldots,\mathbf{e}_d\}$. An example of 1-graded relaxed monomial order is the graded lexicographic order

\begin{proposition}\label{prop:graded-every-d}
    Let $d$ be a positive integer and $\preceq$ be a 1-graded relaxed monomial order in $\mathbb{N}^d$. Assume $\mathbf{e}_1\prec \mathbf{e}_2\prec \cdots \prec \mathbf{e}_d$. Let $S\in \operatorname{R}_\preceq(\mathcal{S}_d)$ and suppose $\operatorname{H}(S)\cap \{\mathbf{e}_1,\ldots,\mathbf{e}_d\}=\{\mathbf{e}_{j_1}\prec \mathbf{e}_{j_2}\prec\ldots\prec\mathbf{e}_{j_r}\}$ for some $\{j_1,j_2,\ldots,j_r\}\subseteq \{1,\ldots,d\}$ and $r\leq d$. Then $\{\mathbf{e}_{j_1}, \mathbf{e}_{j_2},\ldots,\mathbf{e}_{j_r}\}=\{\mathbf{e}_{1},\mathbf{e}_{2},\ldots,\mathbf{e}_r\}$. 
\end{proposition}
\begin{proof}
    If $r=d$ is trivial, so assume $r<d$. By assumptions we have $\operatorname{H}(S)=\{\mathbf{e}_{j_1}\prec \mathbf{e}_{j_2}\prec\ldots\prec\mathbf{e}_{j_r}\prec \cdots\}$. Moreover, by Remark~\ref{rem:trasposizioni}, there exists $\sigma\in \operatorname{P}_d$ such that $\sigma(\{\mathbf{e}_{j_1},\mathbf{e}_{j_2},\ldots,\mathbf{e}_{j_r}\})=\{\mathbf{e}_{1},\mathbf{e}_{2},\ldots,\mathbf{e}_r\}$. In particular, by hypothesis, we have $\operatorname{H}(\sigma(S))=\{\mathbf{e}_{1}\prec \mathbf{e}_{2} \prec \ldots\prec \mathbf{e}_r\prec \cdots\}$. This means that $\sigma(S)\preceq_{\operatorname{R}} S$, and since $S\in \operatorname{R}_\preceq(\mathcal{S}_d)$ we obtain $\sigma(S)=S$. Therefore, $\{\mathbf{e}_{j_1}, \mathbf{e}_{j_2},\ldots,\mathbf{e}_{j_r}\}=\{\mathbf{e}_{1},\mathbf{e}_{2},\ldots,\mathbf{e}_r\}$.
\end{proof}

We find, in the examples that we computed, that the same consequence of Proposition~\ref{prop:graded-every-d} is verified also with some relaxed monomial order not 1-graded. We ask if Proposition~\ref{prop:graded-every-d} holds in general,  without the assumption that $\preceq$ is 1-graded.

\subsection{A different construction for producing $\operatorname{R}_{\preceq}(\mathcal{S}_{g,d})$}
%Now, our aim is to show that it is possible to restrict the transform $\mathcal{O}_\prec$ on $\operatorname{R}_{\preceq}(\mathcal{S}_{g,d})$, for opportune choices of the relaxed monomial order.
In \cite[Section 6]{cisto2021algorithms}, it is shown how to define a rooted tree whose vertices are the GNSs in $\mathcal{S}_{g,d}$ and this construction is related to the transform $\mathcal{O}_\preceq: \mathcal{S}_{g,d} \rightarrow \mathcal{S}_{g,d}$ with  $\mathcal{O}_\preceq(S)=(S\cup \{\mathbf{F}_\preceq(S)\})\setminus \{\mathbf{m}_\preceq(S)\}$. Unfortunately, the transform $\mathcal{O}_\preceq$ cannot be always restricted to $\operatorname{R}_\preceq(\mathcal{S}_{g,d})$. In fact, for $T\in \operatorname{R}_\preceq(S_{g,d})$, the semigroup $T\setminus \{\mathbf{m}_\preceq(S)\}$ could not belong to $\operatorname{R}_\preceq(S_{g,d})$, as shown in the following example.

\begin{example}
Let $T=\mathbb{N}^2\setminus \{(0,1),(1,0),(0,2),(1,1),(2,0),(0,3),(1,2),(3,0),(4,0)\}$ and let $\preceq$ be the graded lexicographic order in $\mathbb{N}^2$. %consider in $\mathbb{N}^2$ the graded lexicographic order $\leq_{\mathrm{GLex}}$. \operatorname{R}_{\leq_{\mathrm{GLex}}}(\mathcal{S}_{9,2})$
Observe that $T\in \operatorname{R}_{\preceq}(\mathcal{S}_{9,2})$, but $\mathbf{m}_\preceq(T)=(2,1)$ and $T\setminus \{(2,1)\}\notin \operatorname{R}_{\preceq}(\mathcal{S}_{10,2})$. Furthermore, if $S=\mathbb{N}^2\setminus \{(0,1),(1,0),(0,2),(1,1),(2,0),(0,3),(1,2),(3,0),(4,0),(4,1)\}\in  \operatorname{R}_{\preceq}(\mathcal{S}_{10,2})$, we have $\mathcal{O}_\preceq(S)=(S\cup \{(4,1)\})\setminus \{(2,1)\}=T\setminus \{(2,1)\}$. Hence, using the graded lexicographic monomial order, it is not possible to restrict the transform $\mathcal{O}_\preceq$ on $\operatorname{R}_{\preceq}(\mathcal{S}_{g,d})$. 
\end{example}

Nevertheless, it is possible to restrict the transform $\mathcal{O}_\preceq$ on $\operatorname{R}_{\preceq}(\mathcal{S}_{g,d})$, for appropriate choices of the relaxed monomial order.

\begin{remark}\label{rem:good-order}
Let $S\in \operatorname{R}_\preceq(\mathcal{S}_{g,d})$ be for positive integers $g,d$. Considering the relaxed monomial orders $\prec_{\text{lex}}$ or $\prec_{1}$ as in Example~\ref{def:orders}, it is not difficult to see that, for all $\sigma\in \operatorname{P}_d$, the inequalities $\mathbf{m}_\preceq(S)\preceq_{\text{lex}}\sigma(\mathbf{m}_\preceq(S))$ and $\mathbf{m}_\preceq(S)\preceq_{1}\sigma(\mathbf{m}_\preceq(S))$ hold. In particular, if $\preceq$ is one of the relaxed monomial orders $\prec_{\text{lex}}$ or $\prec_{1}$, by Lemma~\ref{lem:parent-is-R} and Proposition~\ref{prop:R-remove-gen1} we have $(S\cup \{\mathbf{F}_\preceq(S)\})\setminus \{\mathbf{m}_\preceq (S)\}\in \operatorname{R}_\preceq(\mathcal{S}_{g,d})$.
\end{remark}

%\textcolor{red}{FARE REMARK CHE CON OPPORTUNI ORDINI MONOMIALI LA MOLTEPLICITA' VA SEMPRE BENE}

\medskip
\begin{definition}
Let $g,d$ be positive integers. We say that a relaxed monomial order $\preceq$ in $\mathbb{N}^d$ is $\mathcal{O}$-good if $(S\cup \{\mathbf{F}_\preceq(S)\})\setminus \{\mathbf{m}_\preceq (S)\}\in \operatorname{R}_\preceq(\mathcal{S}_{g,d})$ for all 
$S\in \operatorname{R}_\preceq(\mathcal{S}_{g,d})$. 
\end{definition}

%Let $\preceq$ be a relaxed monomial order in $\mathbb{N}^{d}$ and $\mathbf{s}\in \mathbb{N}^{d}$. The set $S=\{\mathbf{x}\in \mathbb{N}^{d}\mid \mathbf{s}\preceq \mathbf{x}\}\cup \{\mathbf{0}\}$ is a generalized numerical semigroup, and we call it an \emph{ordinary} generalized numerical semigroup with respect to $\preceq$.\\

Let $\preceq$ be a relaxed monomial order and $\{\mathbf{0}=\mathbf{s}_{0}\preceq \mathbf{s}_{1}\preceq \cdots \preceq \mathbf{s}_{g}\}$ the list of the first $g+1$ elements in $\mathbb{N}^{d}$, ordered by $\preceq$. We recall that, as introduced in \cite[Section 6]{cisto2021algorithms}, the set $O_{g,d}(\preceq)=\{\mathbf{x}\in \mathbb{N}^{d}\mid \mathbf{s}_{g}\prec \mathbf{x}\}\cup \{\mathbf{0}\}$ is a GNS and it is called the \emph{ordinary generalized numerical semigroup} in $\mathbb{N}^{d}$ of genus $g$, with respect to $\preceq$.
%\end{definition}

\begin{proposition}
Let $g,d$ be positive integers and $\preceq$ be a relaxed monomial order in $\mathbb{N}^d$. Then $O_{g,d}(\preceq)\in \operatorname{R}_\preceq(\mathcal{S}_{g,d})$.
\end{proposition}
\begin{proof}
Let $\sigma \in \operatorname{P}_d$. Denote $O_{g,d}(\preceq)=\mathbb{N}^d \setminus \{\mathbf{h}_1\prec \mathbf{h}_2 \prec \cdots \mathbf{h}_g\}$, in particular $\mathbf{0}\prec \mathbf{h}_1\prec \mathbf{h}_2 \prec \cdots \mathbf{h}_g$ is the list of the first $g+1$ elements in $\mathbb{N}^{d}$, ordered by $\preceq$. We can assume $\mathbb{N}^d \setminus \sigma(O_{g,d}(\preceq))=\{\sigma(\mathbf{h}_{i_1}) \prec \sigma(\mathbf{h}_{i_2})\prec \cdots \prec \sigma(\mathbf{h}_{i_g})\}$. Let $r=\max \{k\in \{1,\ldots,g\}\mid \mathbf{h}_{j}=\sigma(\mathbf{h}_{i_j})\text{ for all }j\in\{1,\ldots,k\}\}$. If $k=g$ then $O_{g,d}(\preceq)=\sigma(O_{g,d}(\preceq))$ and we can conclude. So, assume $k<g$. Suppose $\sigma(\mathbf{h}_{i_{k+1}})\prec \mathbf{h}_{k+1}$. Then, by the structure of the set $\{\mathbf{h}_1\prec \mathbf{h}_2 \prec \cdots \mathbf{h}_g\}$, we have $\sigma(\mathbf{h}_{i_{k+1}})=\mathbf{h}_j$ for some $j<k+1$. This implies that $\sigma(\mathbf{h}_{i_{k+1}})=\sigma(\mathbf{h}_{i_j})$ with $j\neq k+1$, a contradiction. Therefore, $\mathbf{h}_{{k+1}}\prec \sigma(\mathbf{h}_{i_{k+1}})$. This means that, for every $\sigma\in \operatorname{P}_d$, we have $O_{g,d}(\preceq)\preceq_{\operatorname{R}} \sigma(O_{g,d}(\preceq))$. So, we can conclude.   
\end{proof}

Let $\preceq$ be an $\mathcal{O}$-good relaxed monomial order. We introduce transform $\mathcal{O}^{\operatorname{R}}_{\preceq}:\operatorname{R}_\preceq(\mathcal{S}_{g,d})\setminus \{O_{g,d}(\preceq)\}\rightarrow \operatorname{R}_\preceq(\mathcal{S}_{g,d})$, defined by $\mathcal{O}^{\operatorname{R}}_{\preceq}(S)=(S\cup \{\mathbf{F}_{\preceq}(S)\})\setminus \{\mathbf{m}_\preceq (S)\}$. Let $\mathcal{K}^{\operatorname{R}}_{g,d,\preceq}$ be the oriented graph having set of vertices $\operatorname{R}_\preceq(\mathcal{S}_{g,d})$ and $(S,T)$ is an edge of $\mathcal{K}^{\operatorname{R}}_{g,d,\preceq}$ if $T=\mathcal{O}^{\operatorname{R}}_{\preceq}(S)$. %In such a case we say also that $S$ is a child of $T$.

\begin{theorem}\label{thm:tree2}
Let $g,d\in \mathbb{N}$ and $\preceq$ be an $\mathcal{O}$-good relaxed monomial order. The graph $\mathcal{K}_{g,d,\preceq}^{d}$ is a rooted tree with root $O_{g,d}(\preceq)$. Moreover, if $T\in \mathcal{S}_{g,d}$, then the children of $T$ are the semigroups of the form $T_{\mathbf{h},\mathbf{x}}=(T \cup \{\mathbf{h}\})\setminus \{\mathbf{x}\}$ with $\mathbf{h}\in \operatorname{SG}(T)$, $\mathbf{h}\prec \mathbf{m}_\preceq(T)$, and $\mathbf{x}\in \mathbf{U}_\preceq(T \cup \{\mathbf{h}\})\setminus\{ \mathbf{h}\}$ such that $T_{\mathbf{h},\mathbf{x}}\in \operatorname{R}_\preceq(\mathcal{S}_{g,d})$.
%obtained considering first all $\mathbf{h}\in \operatorname{SG}(T)$ with $\mathbf{h}\preceq \mathbf{m}_{\preceq}(T)$ and after, for all semigroups $T_{\mathbf{h}}=T\cup \{\mathbf{h}\}$, considering $T_{\mathbf{h}}\setminus \{\mathbf{x}\}$ for all $\mathbf{x}\in \operatorname{U}_\preceq(T_{\mathbf{h}})$ and $\mathbf{x}\neq \mathbf{m}_{\preceq}(T_{\mathbf{h}})$.
\label{generagd}
\end{theorem}
\begin{proof}
By the definition of $\mathcal{O}$-good relaxed monomial order, the proof is the same as \cite[Theorem 24]{cisto2021algorithms}.
%Let $T\in \mathcal{S}_{g,d}$, we define the following sequence: 
%\begin{itemize}
%\item $T_{0}=T$,
%\item $T_{i+1}=\left\lbrace\begin{array}{ll}
%         \mathcal{O}_\preceq(T_{i}) & \mbox{if }T_{i}\neq R_{g,d}(\preceq),  \\
%         R_{g,d}(\preceq) & \mbox{otherwise,}
%        \end{array}
%        \right.$
%\end{itemize}
%in particular $T_{i}=\mathcal{O}_\preceq^{i}(T)$ for all $i$. 
%%We know that there exists a nonnegative integer $k$ such that $T_{k}=\mathcal{O}_\preceq^{k}(T)=R_{g,d}(\preceq)$. 
%This sequence stabilizes at a certain point, since the multiplicity must be among the first $g+1$ elements of $\mathbb{N}^d$ (ordered by $\preceq$), and at each step the multiplicity increases. Let $k$ be the first integer such that $T_{k}=R_{g,d}(\preceq)$.
%%
%So the edges $(T_{0},T_{1}),(T_{1},T_{2}),\ldots,(T_{k-1},T_{k})$ provide the unique path from $T$ to the ordinary numerical semigroup of genus $g$ (with respect to $\preceq$).
%
%By Lemma~\ref{treeO2d}, every pair $(T_{\mathbf{h},\mathbf{x}},T)$ is an edge of $\mathcal{T}_{g,\preceq}^{d}$ for every choice of $\mathbf{h}$ and $\mathbf{x}$ as in the statement. Lemma~\ref{treeO1d} ensures that all descendants are of this form.
\end{proof}

\begin{example}\label{exa:tree2}
Let $\preceq$ be the lexicographic order in $\mathbb{N}^{2}$. We compute the tree $\mathcal{K}_{3,2,\preceq}^{\operatorname{R}}$. Consider the semigroup $O_{3,2}(\preceq)=\mathbb{N}^{2}\setminus \{(0,1),(0,2),(0,3)\}$, We have $\mathbf{m}_{\preceq}(O_{3,2}(\preceq))=(0,4)$, so $\{\textbf{h}\in\operatorname{SG}(O_{3,2}(\preceq))\mid \textbf{h}\preceq \mathbf{m}_{\preceq}(0,4))\}=\{(0,2),(0,3)\}$. Consider the semigroups:
\begin{itemize}
\item $T_{(0,2)}=O_{3,2}(\preceq)\cup \{(0,2)\}=\mathbb{N}^{2}\setminus \{(0,1),(0,3)\}$, with $\operatorname{U}_\preceq(T_{(0,2)})=\{(1,0), (1,1),(0,5)\}$,  
\item $T_{(0,3)}=O_{3,2}(\preceq)\cup \{(0,3)\}=\mathbb{N}^{2}\setminus \{(0,1),(0,2)\}$, with $\operatorname{U}_\preceq(T_{(0,3)})=\{(1,0),(1,1),(0,3),(1,2),(0,4),(0,5)\}$  
\end{itemize}

So we can find the children of $O_{3,2}(\preceq)$ among these semigroups:
\begin{itemize}
\item $S_{1}=T_{(0,2)}\setminus \{(1,0)\}=\mathbb{N}^{2}\setminus \{(0,1),(0,3),(1,0)\}\in \operatorname{R}_{\preceq}(\mathcal{S}_{g,d})$,
\item $S_{2}=T_{(0,2)}\setminus \{(1,1)\}=\mathbb{N}^{2}\setminus \{(0,1),(0,3),(1,1)\} \in \operatorname{R}_{\preceq}(\mathcal{S}_{g,d})$,
\item $S_{3}=T_{(0,2)}\setminus \{(0,5)\}=\mathbb{N}^{2}\setminus \{(0,1),(0,3),(0,5)\} \in \operatorname{R}_{\preceq}(\mathcal{S}_{g,d})$,
\item $S_{4}=T_{(0,3)}\setminus \{(1,0)\}=\mathbb{N}^{2}\setminus \{(0,1),(0,2),(1,0)\}\in \operatorname{R}_{\preceq}(\mathcal{S}_{g,d})$,
\item $S_{5}=T_{(0,3)}\setminus \{(1,1)\}=\mathbb{N}^{2}\setminus \{(0,1),(0,2),(1,1)\}\in \operatorname{R}_{\preceq}(\mathcal{S}_{g,d})$,
\item $S_{6}=T_{(0,3)}\setminus \{(1,2)\}=\mathbb{N}^{2}\setminus \{(0,1),(0,2),(1,2)\}\in \operatorname{R}_{\preceq}(\mathcal{S}_{g,d})$,
\item $S_{7}=T_{(0,3)}\setminus \{(0,4)\}=\mathbb{N}^{2}\setminus \{(0,1),(0,2),(0,4)\}\in \operatorname{R}_{\preceq}(\mathcal{S}_{g,d})$,
\item $S_{8}=T_{(0,2)}\setminus \{(0,5)\}=\mathbb{N}^{2}\setminus \{(0,1),(0,2),(0,5)\}\in \operatorname{R}_{\preceq}(\mathcal{S}_{g,d})$.
\end{itemize}
For $S_{1},S_{2},S_{3},S_{6},S_{7},S_{8}$ the set $\{\mathbf{h}\in \operatorname{SG}(S)\mid \mathbf{h}\preceq \mathbf{m}_{\preceq}(S)\}$ is empty, so they are leaves. Therefore, we need to consider the semigroups $S_{4}$ and $S_{5}$. From $S_{5}$, the only semigroup of kind $T_{\mathbf{h},\mathbf{x}}$ as in Theorem~\ref{thm:tree2} is
 $S_{14}=S_{5}\cup \{(0,2)\}\setminus \{(2,1)\}=\mathbb{N}^{2}\setminus \{(0,1),(1,1),(2,1)\}\notin \operatorname{R}_{\preceq}(\mathcal{S}_{g,d})$. So, $S_5$ is a leaf. From $S_{4}$, the semigroups $T_{\mathbf{h},\mathbf{x}}$ as in Theorem~\ref{thm:tree2} are:
\begin{itemize}
\item $S_{9}=S_{4}\cup \{(0,2)\}\setminus \{(1,1)\}=\mathbb{N}^{2}\setminus \{(0,1),(1,0),(1,1)\} \in \operatorname{R}_{\preceq}(\mathcal{S}_{g,d})$,
\item $S_{10}=S_{4}\cup \{(0,2)\}\setminus \{(2,1)\}=\mathbb{N}^{2}\setminus \{(0,1),(1,0),(2,1)\}\notin \operatorname{R}_{\preceq}(\mathcal{S}_{g,d})$,
\item $S_{11}=S_{4}\cup \{(0,2)\}\setminus \{(1,2)\}=\mathbb{N}^{2}\setminus \{(0,1),(1,0),(1,2)\}\in \operatorname{R}_{\preceq}(\mathcal{S}_{g,d})$,
\item $S_{12}=S_{4}\cup \{(0,2)\}\setminus \{(2,0)\}=\mathbb{N}^{2}\setminus \{(0,1),(1,0),(2,0)\}\notin \operatorname{R}_{\preceq}(\mathcal{S}_{g,d})$,
\item $S_{13}=S_{4}\cup \{(0,2)\}\setminus \{(1,1)\}=\mathbb{N}^{2}\setminus \{(0,1),(1,0),(3,0)\}\notin \operatorname{R}_{\preceq}(\mathcal{S}_{g,d})$.
\end{itemize}
Hence, we have to continue the procedure from the semigroups $S_9$ and $S_{11}$. It is not difficult to verify that $S_11$ is a leaf. From the semigroup $S_{9}$ we obtain:
\begin{itemize}
\item $S_{15}=S_{9}\cup \{(0,1)\}\setminus \{(2,0)\}=\mathbb{N}^{2}\setminus \{(1,0),(1,1),(2,0)\}\notin \operatorname{R}_{\preceq}(\mathcal{S}_{g,d})$,
\item $S_{16}=S_{9}\cup \{(0,1)\}\setminus \{(3,0)\}=\mathbb{N}^{2}\setminus \{(1,0),(1,1),(3,0)\}\notin \operatorname{R}_{\preceq}(\mathcal{S}_{g,d})$,
\item $S_{17}=S_{9}\cup \{(0,1)\}\setminus \{(1,2)\}=\mathbb{N}^{2}\setminus \{(1,0),(1,1),(1,2)\}\in \operatorname{R}_{\preceq}(\mathcal{S}_{g,d})$.
\end{itemize}

The shape of the tree $\mathcal{K}_{3,2,\preceq}^{\operatorname{R}}$ is depicted in Figure~\ref{fig:tree-g-d}. For a comparison with the analogous tree which produces the whole set $\mathcal{S}_{g,d}$, see \cite[Example 26]{cisto2021algorithms}

\begin{figure}[h!]
%    \subfloat[]{\begin{tikzpicture}
% \tikzset{level distance=6em}
% \Tree
%         [.$R_{3,2}(\preceq)$  
%              $S_1$ 
%           $S_2$ 
%                $S_3$
%            [.$S_4$ 
%                       [.$S_9$ $S_{15}$ $S_{16}$ $S_{17}$   ] 
%                     $S_{10}$
%                      $S_{11}$
%                       [.$S_{12}$ $S_{18}$ $S_{19}$ $S_{20}$ $S_{21}$ ]
%                  [.$S_{13}$ $S_{22}$ ]
%          ] 
%           [.$S_5$ $S_{14}$ ]               
%                $S_6$
%            $S_7$
%                $S_8$
%    ]
% \end{tikzpicture} } 
\begin{tikzpicture}
\tikzset{level distance=6em}
\Tree
        [.$R_{3,2}(\preceq)$  
             $S_1$ 
          $S_2$ 
               $S_3$
           [.$S_4$ 
                      [.$S_9$ $S_{17}$   ] 
                     $S_{11}$
         ] $S_5$ $S_{14}$                
               $S_6$
           $S_7$
               $S_8$
   ]
\end{tikzpicture}
    \caption{The tree $\mathcal{K}_{3,2,\preceq}^{\operatorname{R}}$, where $\preceq$ is the lexicographic order, with reference to Example~\ref{exa:tree2}}
    \label{fig:tree-g-d}
\end{figure}
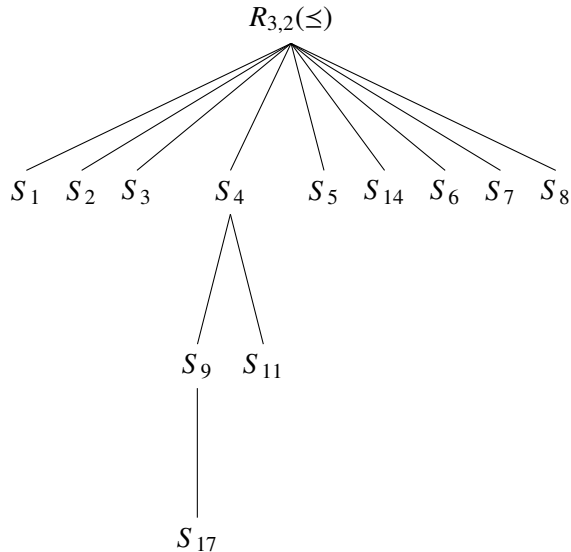
\end{example}

\begin{algorithm}[H]  \label{alg1}
\caption{Another algorithm for computing the set $\operatorname{R}_\preceq(\mathcal{S}_{g,d})$ \label{alg:S-g-d}}
\DontPrintSemicolon

\KwData{Two integers $g,d\in \mathbb{N}$ and a relaxed monomial order $\preceq$ in \(\mathbb{N}^d\).}
\KwResult{$\operatorname{R}_\preceq(\mathcal{S}_{g,d})$}

 Compute $O_{g,d}(\preceq)$.\; 
 $\mathcal{P}=\{O_{g,d}(\preceq)\}, \mathcal{L}=\{O_{g,d}(\preceq)\}$.\;
\nl \While{There exists $S\in \mathcal{L}$ \emph{such that} \(\{\mathbf{h}\in \operatorname{SG}(S)\mid \mathbf{h}\preceq \mathbf{m}_{\preceq}(S)\}\neq \emptyset\)}{
   $\mathcal{I}=\emptyset$\;
   \For{$S\in \mathcal{L}$}
                {
       \If{\(\{\mathbf{h}\in \operatorname{SG}(S)\mid \mathbf{h}\preceq \mathbf{m}_{\preceq}(S)\}\neq \emptyset\)}{
          $\mathcal{R}=\emptyset$\;
         \For{$\mathbf{h}\in \{\mathbf{h}\in \operatorname{SG}(S)\mid \mathbf{h}\preceq \mathbf{m}_{\preceq}(S)\}$}{
             $\mathcal{R}=\mathcal{R}\cup \{S\cup\{\mathbf{h}\}\}$
           } 
         \nl    \For{$T\in \mathcal{R}$}{
                  \For{$\mathbf{x}\in \mathbf{U}_\preceq(T)$ with $\mathbf{x}\neq \mathbf{m}_{\preceq}(S)$}{
                  \nl \If{$T\setminus \{\mathbf{x}\}=\min_{\preceq_{\operatorname{R}}}([T\setminus \{\mathbf{x}\}]_\simeq)$}{
                  $\mathcal{I}=\mathcal{I}\cup \{T\setminus\{\mathbf{x}\}\}$\;}
                  %$\mathcal{S}_{g,d}=\mathcal{S}_{g,d}\cup L$
                      }
                      }
                   }
                 }
                 $\mathcal{P}=\mathcal{P}\cup \mathcal{I}$\;
                 $\mathcal{L}=\mathcal{I}$\;
                }
 \Return $\mathcal{P}$\;  
\end{algorithm}

\medskip
Some questions about the nature of $\mathcal{O}$-good orders naturally arises. We do not know if it is possible to characterize what are the $\mathcal{O}$-good relaxed monomial order, in terms of some particular properties. In our examples, we obtained that in each $\mathcal{O}$-good order $\preceq$ it is verified $\mathbf{m}_\preceq(S)\preceq \sigma(\mathbf{m}_\preceq(S))$ for every $S\in \mathcal{S}_{g,d}$ and for every $\sigma \in \operatorname{P}_d$. We ask if there exists an $\mathcal{O}$-good order such $\sigma(\mathbf{m}_\preceq(S))\prec \mathbf{m}_\preceq(S)$ for some $\sigma\in \operatorname{P}_d$ and $S\in \operatorname{R}_\preceq(\mathcal{S}_{g,d})$.

\section{On the number of GNSs up to isomorphism and computational data}

%Let $g,d$ be positive integers. Denote $N_{g,d}=|\mathcal{S}_{g,d}|$ and  $N'_{g,d}=|\operatorname{R}_\preceq(\mathcal{S}_{g,d})|$.

Let $A$ be a subset of $\mathbb{N}^{d}$, we denote by $\mathrm{Span}_\mathbb{R}(A)$ the $\mathbb{R}$-vector subspace of $\mathbb{R}^{d}$ generated by the elements of $A$. Recall that a vector subspace of $\mathbb{R}^{d}$ is a \emph{coordinate linear space} if it is spanned by a subset of the set of standard basis vectors of $\mathbb{R}^d$. %$\{\textbf{e}_{1},\textbf{e}_{2},\ldots,\textbf{e}_{d}\}$.  

\begin{proposition}[\cite{failla2016algorithms}, Proposition 5.2] Let $S\subseteq \mathbb{N}^{d}$ be a generalized numerical semigroup and $\operatorname{H}(S)$ the set of its holes. Then $\mathrm{Span}_{\mathbb{R}}(\operatorname{H}(S))$ is a coordinate linear space.\label{spaziocoord}\end{proposition}

Let $g,d$ be positive integers, consider the following notations:
\begin{itemize}
%\item $S_{g,d}$ is the set of all generalized numerical semigroup with genus $g$ in $\mathbb{N}^{d}$%
%\item $\mathcal{S}_{g,d}^{(r)}=\{S\in S_{g,d}\ |\ \dim(\mathrm{Span}_{\mathbb{R}}(\operatorname{H}(S)))=r\}$.

%\item $\mathcal{R}_{g,d}^{(r)}=\{[S]_\simeq \mid S\in \mathcal{S}_{g,d}^{(r)}\}$.

\item $N_{g,d}=|\{[S]_\simeq \mid S\in S_{g,d}\}|$ and $n_{g,d}=|\mathcal{S}_{g,d}|$.

\item $\mathcal{R}_{g,d}^{(r)}=\{[S]_\simeq \mid S\in S_{g,d}\text{ and }\dim(\mathrm{Span}_{\mathbb{R}}(\operatorname{H}(S)))=r\}$ and $N^{(r)}_{g,d}=\lvert \mathcal{R}_{g,d}^{(r)}\rvert$, with $r\leq \min\{g,d\}$.

%\item $\operatorname{R}_\preceq \left (\mathcal{S}_{g,d}^{(r)}\right)=\{S\in \operatorname{R}_\preceq(S_{g,d})\ |\ \dim(\mathrm{Span}_{\mathbb{R}}(\operatorname{H}(S)))=r\}$.
%\item $N_{g,d}$ and $N_{g,d}^{(r)}$ denote respectively the cardinalities of $S_{g,d}$ and $S_{g,d}^{(r)}$.

%\item $N_{g,d}=|\operatorname{R}_\preceq(\mathcal{S}_{g,d})|$ and $N^{(r)}_{g,d}=|\operatorname{R}_\preceq(\mathcal{S}_{g,d})|$
\end{itemize}

%Let $g,d$ be positive integers. Denote $N_{g,d}=|\mathcal{S}_{g,d}|$ and  $N'_{g,d}=|\operatorname{R}_\preceq(\mathcal{S}_{g,d})|$.

Observe that $N_{g,d}=|\operatorname{R}_\preceq(\mathcal{S}_{g,d})|$ for every relaxed monomial order $\preceq$.

\begin{proposition}\label{prop:sum_Ngi}
Let $g,d$ be positive integers and $q=\min \{g,d\}$. Then $$ N_{g,d}=\sum_{n=1}^q N^{(n)}_{g,n}.$$

\end{proposition}
\begin{proof}
First of all, observe that $\{[S]_\simeq \mid S\in S_{g,d}\}=\bigsqcup_{n=1}^q \mathcal{R}_{g,d}^{(n)}$. So, in order to obtain our results, it is sufficient to prove that $|\mathcal{R}_{g,d}^{(n)}|=|\mathcal{R}_{g,n}^{(n)}|$ for all $n\in \{1,\ldots,q\}$. Fix $n\in \{1,\ldots,q\}$. If $n=d$ the claim is trivial, so we can assume $n<d$. Let us denote with $\{\mathbf{e}_1,\ldots,\mathbf{e}_d\}$ and $\{\mathbf{e}'_1,\ldots,\mathbf{e}'_n\}$ the sets of the standard basis vectors in $\mathbb{R}^d$ and $\mathbb{R}^n$, respectively. If $S\in \mathcal{R}_{g,d}^{(n)}$, then there exists $\{i_1,\ldots,i_n\}\subseteq \{1,\ldots,d\}$ such that for all $\mathbf{h}=\sum_{j=1}^d h_j\mathbf{e}_j\in \operatorname{H}(S)$ we have $h_k=0$ for all $k\notin \{i_1,\ldots,i_n\}$. In particular, if $\mathbf{h}\in \operatorname{H}(S)$, then we can express $\mathbf{h}=\sum_{j=1}^n h_{i_j}\mathbf{e}_{i_j}$ and we define $\mathbf{h}'= \sum_{j=1}^n h_{i_j}\mathbf{e}'_{j}\in \mathbb{N}^n$. Now, we denote  $\overline{\operatorname{H}}(S)=\{\mathbf{h}'\in \mathbb{N}^d\mid \mathbf{h}\in \operatorname{H}(S)\}$. It is not difficult to show that the set $\overline{S}=\mathbb{N}^n \setminus \overline{\operatorname{H}}(S)$ is a GNS. In particular, $\operatorname{H}(\overline{S})=\overline{\operatorname{H}}(S)$. Having in mind these notations, we define the following:

$$ \psi : \mathcal{R}_{g,d}^{(n)} \longrightarrow \mathcal{R}_{g,n}^{(n)}\qquad \text{ with }\qquad [S]_\simeq \longmapsto [\overline{S}]_\simeq $$

\noindent Firstly, we show that $\psi$ is well defined. Let $S_1,S_2$ be two GNSs such that $[S_1]_\simeq = [S_2]_\simeq$ in $\mathcal{R}^{(n)}_{g,d}$. This means that there exists $\sigma\in \operatorname{P}_d$ such that $\sigma(S_2)=S_1$. By Proposition~\ref{prop:iso-gaps}, we have also $\sigma(\operatorname{H}(S_2))=\operatorname{H}(S_1)$. So, assuming $\mathrm{Span}_{\mathbb{R}}(\operatorname{H}(S_2))= \mathrm{Span}_{\mathbb{R}}(\{\mathbf{e}_{i_1},\ldots,\mathbf{e}_{i_n}\})$ with $\{i_1,\ldots,i_n\}\subseteq \{1,\ldots,d\}$, we have $\mathrm{Span}_{\mathbb{R}}(\operatorname{H}(S_1))= \mathrm{Span}_{\mathbb{R}}(\{\mathbf{e}_{\sigma(i_1)},\ldots,\mathbf{e}_{\sigma(i_n)}\})$. In particular, considering the permutation $\overline{\sigma}\in \operatorname{P}_n$ defined by $\overline{\sigma}(k)=\ell$ such that $\sigma(i_k)=i_\ell$ for all $k\in\{1,\ldots,n\}$, it is not difficult to see that $\overline{\sigma}(\operatorname{H}(\overline{S}_2))=\operatorname{H}(\overline{S}_1)$. Hence, $\overline{\sigma}(\overline{S}_2)=\overline{S}_1$, obtaining $[\overline{S}_1]_\simeq = [\overline{S}_2]_\simeq$. Therefore $\psi$ is a function. Now, we show that $\psi$ is injective. Let $S_1,S_2\in \mathcal{S}_{g,d}$ such  that $\mathrm{Span}_{\mathbb{R}}(\operatorname{H}(S_1))= \mathrm{Span}_{\mathbb{R}}(\{\mathbf{e}_{i_1},\ldots,\mathbf{e}_{i_n}\})$ and $\mathrm{Span}_{\mathbb{R}}(\operatorname{H}(S_2))= \mathrm{Span}_{\mathbb{R}}(\{\mathbf{e}_{j_1},\ldots,\mathbf{e}_{j_n}\})$, with $\{i_1,\ldots,i_n\}$ and $\{j_1,\ldots,j_n\}$ subsets of $\{1,\ldots,d\}$. Assume that $[\overline{S}_1]_\simeq =[\overline{S}_2]_\simeq$. So, there exists $\overline{\sigma}\in \operatorname{P}_n$ such that $\overline{\sigma}(\operatorname{H}(\overline{S}_2))=\operatorname{H}(\overline{S}_1)$. As a consequence, we have $\{i_1,\ldots,i_n\}=\{j_1,\ldots,j_n\}$. In particular, if $\sigma$ is any partition on $\{1,\ldots,d\}$ such that $\sigma(i_k)=j_{\overline{\sigma}(k)}$ for all $k\in \{1,\ldots,n\}$, then $\sigma(\operatorname{H}(S_2))=\operatorname{H}(S_1)$. This means that  $[S_1]_\simeq = [S_2]_\simeq$, so $\psi$ is injective. Finally, we show that $\psi$ is surjective. Let $T'\in \mathcal{R}_{g,n}^{(n)}$. For all $\mathbf{h}'=\sum_{j=1}^n h_j' \mathbf{e}'_j\in \operatorname{H}(T')$ we define $\mathbf{h}=\sum_{j=1}^d h_j \mathbf{e}_j\in \mathbb{N}^d$, where $h_j=h'_j$ for all $j\in \{1,\ldots,n\}$ and $h_j=0$ for all $j\in \{n+1,\ldots,d\}$. Following this definition, define also $T=\mathbb{N}^d\setminus \{\mathbf{h}\in \mathbb{N}^d\mid \mathbf{h}'\in \operatorname{H}(T')\}$. It is not difficult to see that $T\subseteq \mathbb{N}^d$ is a GNS such that $T\in \mathcal{R}^{(n)}_{g,d}$ and $\psi([T]_\simeq)=[T']_\simeq$. So, $\psi$ is surjective. In conclusion, we have that $\psi$ is bijective. As a consenquence, $|\mathcal{R}_{g,d}^{(n)}|=|\mathcal{R}_{g,n}^{(n)}|$ and this completes the proof.
\end{proof}

\begin{corollary}
Let $g$ be a positive integer. Then $N_{g,d}=N_{g,g}$ for all $d\geq g$.
\end{corollary}
\begin{proof}
If $d\geq g$ then $q=\min \{d,g\}=g$. So, by Proposition~\ref{prop:sum_Ngi}, $ N_{g,d}=\sum_{n=1}^g N^{(n)}_{g,n}$. But it holds trivially $q=g=\min\{g,g\}$. In particular, Proposition~\ref{prop:sum_Ngi} gives also $N_{g,g}=\sum_{n=1}^g N^{(n)}_{g,n}$. Therefore, $N_{g,d}=N_{g,g}$.
\end{proof}

% \begin{table}
% \label{tabella1}
% \caption{Computational results for $N_{g,2}$}
% \begin{tabular}{|c|cc|cc|}
% \toprule
% $g$ &  $n_{g,2}$ & $N_{g,2}$ & $n_{g,3}$ & $N_{g,3}$ \\
% \midrule

% 1	& 2	&1	&  3 &	1  \\

% 2 &	7	& 4	& 15	& 4 \\
% 3 &	23	& 12  &	67	&  15\\
% 4 &	71	& 37 &	292  & 	59\\
% 5 &	210	 & 107  &  1215  &	224\\
% 6 &	638	 & 323 &  5075  &	903\\
% 7 &	1894 &	953  &	20936  &	3611\\
% 8 &	5570 &	2798 &	85842  &	14603\\
% 9 &	16220 &	8128 &	349731 &	58954\\
% 10 & 46898 & 23486 & 1418323 &	237956\\
% 11 & 134856 & 67477 &	5731710 &	\\
% 12 & 386354 &	193285 & 23100916 &	\\
% 13 & 1102980 &	551628 & 92882954 &	\\
% 14 & 3137592 &	1569107 & 372648740 &	\\
% %15 & 8892740 &	4240572 &	0,4768577514 &	2,8342563342\\
% %16	& 25114649	& 12030332	& 0,4790165294	& 2,8241744389\\
% %17	& 70686370	& 34007389	& 0,4811024954	& 2,8145473982\\
% %18 & 198319427 & 95801019	& 0,4830642184	& 2,8056247194\\
% %19 & 554813870	& 269005797	 & 0,4848577362	 & 2,797577012\\
% %20	& 1548231268 & 753133297 & 0,4864475434	& 2,7905417505\\
% %21	& 4310814033 & 2103045138	& 0,4878533664	& 2,7843476114\\

% \bottomrule
% \end{tabular}
% \qquad
% \begin{tabular}{|c|ccc|}
% \toprule
% $g$ &  $N_{g,4}$ & $N_{g,5}$ & $N_{g,6}$  \\
% \midrule

%  1 & 1 &  1  &	   1 \\
% 2 &  4 &  4 &	  4 \\
% 3 &  15 & 15 &	 15 \\	
% 4 & 64 & 64 &    64 \\
% 5 & 270 &	277&	277  \\
% 6 & 1254 &	1344&	1355 \\
% 7 & 5945 &	6810&	\\
% 8 & 29132 &	36536 & \\

% \bottomrule
% \end{tabular}

% \end{table}

Some values of $N_{g,d}$ are provided in Table~\ref{tabella1} by implementation of Algorithm~\ref{alg:ngd} in the computer algebra software \texttt{GAP}\cite{GAP} and with the help of the package \texttt{numericalsgps} \cite{numericalsgps}.

\begin{table}[h!]
\begin{tabular}{|c|cc|cc|ccc|}
\toprule
$g$ &  $n_{g,2}$ & $N_{g,2}$ & $n_{g,3}$ & $N_{g,3}$ & $N_{g,4}$ & $N_{g,5}$ & $N_{g,6}$ \\
\midrule

1	& 2	&1	&  3 &	1 & 1 &  1  &	   1 \\

2 &	7	& 4	& 15	& 4 &  4 &  4 &	  4 \\
3 &	23	& 12  &	67	&  15 &  15 & 15 &	 15\\
4 &	71	& 37 &	292  & 	59 & 64 & 64 &    64\\
5 &	210	 & 107  &  1215  &	224 & 270 &	277&	277\\
6 &	638	 & 323 &  5075  &	903 & 1254 &	1344&	1355\\
7 &	1894 &	953  &	20936  &	3611 & 5945 &	6810&\\
8 &	5570 &	2798 &	85842  &	14603 & 29132 &	36536 &\\
9 &	16220 &	8128 &	349731 &	58954 & & &\\
10 & 46898 & 23486 & 1418323 &	237956 & & &\\
11 & 134856 & 67477 &	5731710 &	& & &\\
12 & 386354 &	193285 & 23100916 &	& & &\\
13 & 1102980 &	551628 & 92882954 &	& & &\\
14 & 3137592 &	1569107 & 372648740 &	& & &\\
\bottomrule
\end{tabular}
\caption{Some computational data of $N_{g,d}$ (compared with $n_{g,d}$ for $d=2,3$)}
\label{tabella1}
\end{table}

Finally, we conclude with some possible developments that arises from this work. In the context of numerical semigroups and affine semigroups there are different graph-tree constructions to produce families of semigroups with prescribed properties. For instance, we mention the following in the context of affine semigroups:

\begin{itemize}
    \item In \cite{bernardini_corner}, the authors introduce the notion of \emph{corner element} in the context of GNSs, as a generalization of the concept of \emph{conductor} of a numerical semigroup. The authors also provide a procedure to generate all GNSs with a given corner element.

    \item In \cite{garcia_someproperties, garcia2024ideals}, the authors consider the notion of $\mathcal{C}$-semigroup, that is a an affine semigroup having finite complement in a integer cone $\mathcal{C}$. A GNS can be seen as a $\mathcal{C}$-semigroup with $\mathcal{C}=\mathbb{N}^d$. In these articles it is shown, respectively, a procedure to compute all irreducible $\mathcal{C}$-semigroups of given Frobenius element and how to arrange the set of the so called \emph{I(S)-semigroups} in a rooted tree.

    \item Another construction is given in \cite{cisto2021almost}, in order to produce all almost symmetric GNSs with a given Frobenius element.
\end{itemize}

A question that arises, is if it is possible to refine these procedures in order to generated the prescribed families of semigroups, in such a way that only one semigroup in each isomorphism equivalence class is produced.

% \medskip 
% \textcolor{red}{
% \textbf{Questions}: 
% \begin{itemize}
% \item For $T\in \operatorname{R}_\preceq(\mathcal{S}_{g,d})$ and $\mathbf{n}$ is a minimal generator of $T$, is there a characterization, ``easy to check (computationally)'', in order to $T\setminus \{\mathbf{n}\}\in \operatorname{R}_\preceq(\mathcal{S}_{g,d})$? Is there a characterization at least for generators bigger than the Frobenius element with respect to $\preceq$?
% \item We ask if there exists an $\mathcal{O}$-good order such $\sigma(\mathbf{m}_\preceq(S))\prec \mathbf{m}_\preceq(S)$ for some $\sigma\in \operatorname{P}_d$ and $S\in \operatorname{R}_\preceq(\mathcal{S}_{g,d})$ .
% \item Characterize $\mathcal{O}$-good relaxed monomial order. For instance, is there any non locally finite relaxed monomial order which is an $\mathcal{O}$-good order?
% \end{itemize}}

% \begin{table}[h!]
% \label{tabella2}
% \caption{Some values of $N_{g,d}$ for $d>2$}
% \begin{tabular}{|c|ccc|}
% \toprule
% $g$ &  $N_{g,4}$ & $N_{g,5}$ & $N_{g,6}$  \\
% \midrule

%  1 & 1 &  1  &	   1 \\
% 2 &  4 &  4 &	  4 \\
% 3 &  15 & 15 &	 15 \\	
% 4 & 64 & 64 &    64 \\
% 5 & 270 &	277&	277  \\
% 6 & 1254 &	1344&	1355 \\
% 7 & 5945 &	6810&	\\
% 8 & 29132 &	36536 & \\

% \bottomrule
% \end{tabular}
% \end{table}
%\preccurlyeq
\subsection*{Acknowledgements} All authors acknowledge support of the group GNSAGA of INdAM. The third author is supported by Scientific and Technological Research Council of Turkey T\"UB\.{I}TAK under the Grant No: 122F128, and he is thankful to T\"UB\.{I}TAK for their supports.

%
%\nocite{*}
\bibliographystyle{plain}
\bibliography{miabiblio}

%\begin{thebibliography}{20}

%\end{thebibliography}

\end{document}